\documentclass[11pt]{amsart}
\usepackage[margin=1.2in]{geometry}
\usepackage[utf8]{inputenc}
\usepackage{graphicx}
\usepackage{amssymb, amsmath, amsthm, graphics}
\usepackage{mathabx,epsfig}
\usepackage[mathscr]{euscript}
\usepackage{tikz}
\usetikzlibrary{calc}
\usetikzlibrary{matrix}
\usepackage{latexsym}
\usepackage[all]{xy}
\usepackage{color}

\usepackage{bbold}

\usepackage{tikz}
\usepackage{tikz-cd}
\input xy
\xyoption{all}
\usepackage{euscript}
\usepackage{multirow}
\usepackage{etex, pictexwd,dcpic}
\usetikzlibrary{positioning}
\usetikzlibrary{shapes.geometric}
\usetikzlibrary{shapes.misc}
\usetikzlibrary{calc}
\usetikzlibrary{positioning}

\usepackage{pgflibraryarrows}
\usepackage{pgflibrarysnakes}

\usetikzlibrary{trees} % LATEX and plain TEX 
\usetikzlibrary[trees] % ConTEXt

\theoremstyle{plain}
\newtheorem{theorem}{Theorem}[section]
\newtheorem{corollary}[theorem]{Corollary}
\newtheorem{proposition}[theorem]{Proposition}

\newtheorem{lemma}[theorem]{Lemma}
\theoremstyle{definition}
\newtheorem{definition}[theorem]{Definition}
\newtheorem{notation}[theorem]{Notation}

\theoremstyle{remark}
\newtheorem{remark}[theorem]{Remark}

\newtheorem{example}[theorem]{Example}

\numberwithin{equation}{section}

\renewcommand{\(}{\begin{equation*}}
\renewcommand{\)}{\end{equation*}}
\newcommand{\bea}{\begin{eqnarray*}}
\newcommand{\eea}{\end{eqnarray*}}

\def\proof {{Proof.}\hspace{7pt}}

\def\endofproof {\hfill{$\Box$}\\}

\newcommand{\beq}{\begin{equation}}
\newcommand{\eeq}{\end{equation}}

\newcommand{\onto}{\twoheadrightarrow}

\newcommand{\into}{\hookrightarrow}

\numberwithin{equation}{section}

\renewcommand{\(}{\begin{equation}}
\renewcommand{\)}{\end{equation}}

\def\1{{\bf 1}}

\def\<{\langle}
\def\>{\rangle}

\numberwithin{equation}{section}

\newcommand{\RR}{\ensuremath{\mathbb R}}
\newcommand{\NN}{\ensuremath{\mathbb N}}

\newcommand{\ZZ}{\ensuremath{\mathbb Z}}

\newcommand{\CC}{\ensuremath{\mathbb C}}

\newcommand{\sset}{\ensuremath{s\mathscr{S}\mathrm{et}}}

\newcommand{\set}{\ensuremath{\mathscr{S}\mathrm{et}}}

\newcommand{\map}{\mathrm{Map}}
\makeatletter
\newcommand{\colim@}[2]{%
  \vtop{\m@th\ialign{##\cr
    \hfil$#1\operator@font colim$\hfil\cr
    \noalign{\nointerlineskip\kern1.5\ex@}#2\cr
    \noalign{\nointerlineskip\kern-\ex@}\cr}}%
}
\newcommand{\colim}{%
  \mathop{\mathpalette\colim@{\rightarrowfill@\textstyle}}\nmlimits@
}

\makeatother

\def\sm{{\rm C}^{\infty}}

\def\Cart{{\sf Cart}}

\def\Sp{\mathscr{S}{\rm p}}

    \def\Ab{\mathscr{A}{\rm b}}

\title{On the differential $K$-theory of moduli stacks}
\author{Daniel Grady}

\begin{document}

\maketitle

\begin{abstract}
  We compute the connective differential $K$-theory and the differential cohomology of the moduli stack of principal $G$-bundles with connection. The results are formulated in terms of invariant polynomials and the representation ring of $G$. We use the homotopy theory of presheaves of spaces and presheaves of spectra to establish the results. 
\end{abstract}

\setcounter{tocdepth}{1}
\tableofcontents

\section{Introduction}
%Differential cohomology has been a subject of interest since its conception. The theory includes differential form representatives for cocycles at the level of cohomology. This is important in the construction of secondary invariants, such as the Chern--Simons invariant, which was in fact the original motivation for its construction in \cite{CS}. Since its construction, the theory has found many applications to physics. Notably, it was discussions with Edward Witten on quantization of the M5 brane that led Hopkins and Singer to generalize the construction to other cohomology theories in \cite{HS}. 
%
%Despite being of considerable interest in physics, very few computations of differential cohomology groups exist in the literature. Part of the reason for this is that the groups are almost never countably generated. Even over $\RR$, they are almost always infinite dimensional. Although this is the case for generic smooth manifolds, it is often not the case for moduli stacks, since in this case the differential form data must satisfy the strong requirement that it is natural with respect to pullback of whatever objects are being classified by the stack. 

In this note, we compute the differential cohomology (\cite{CS},\cite{BB},\cite{Bunke},\cite{HS},\cite{BNV},\cite{GradySati}) and the connective differential $K$-theory of the moduli stack of smooth principal $G$-bundles with connections $\mathbf{B}_{\nabla}G$. The differential cohomology of $\mathbf{B}_{\nabla}G$ is almost certainly known to experts, although we have not been able to find a precise statement in the literature. Here, we generalize the calculation to differential cohomology with coefficients in an arbitrary graded ring $A$ of finite type, which serves to fill this gap and gives a precise statement in full generality.

Throughout the paper, we let $G$ be an arbitrary Lie group and we denote the symmetric algebra of Ad-invariant polynomials in $\mathfrak{g}^{\vee}$ by $I^*(\mathfrak{g}^{\vee})={\rm Sym}^*(\mathfrak{g}^{\vee})^G$, where the grading is given by twice the degree. Let $A$ be a graded ring of finite type (i.e., finitely generated in each degree). By convention, we take graded ring to mean graded $\ZZ$-algebra, so that $A_0=\ZZ$. Let $I_{A}^*(\mathfrak{g}^{\vee})$ be the $\ZZ$-graded ring of invariant polynomials with coefficients in $A$ (Defintion \ref{invcoA}). 
\begin{theorem}\label{theoremA}
Let $A$ be a graded ring of finite type that is torsion free in each degree (as a $\ZZ$-module). Let $G$ be any Lie group and let $I_{A}^*(\mathfrak{g}^{\vee})_H\subset I_{A}^*(\mathfrak{g}^{\vee})$ denote the subring of $H$-integral elements. Then the curvature map defines a surjective homomorphism of graded rings
$$\mathcal{R}:\widehat{H}^*(\mathbf{B}_{\nabla}G;A)\to I_A^*(\mathfrak{g}^{\vee})_{H}.$$
Moreover, if $G$ is compact and connected, then the forgetful map $\mathcal{I}:\widehat{H}^*(\mathbf{B}_{\nabla}G;A)\to H^*(BG;A)$ maps the kernel of $\mathcal{R}$ isomorphically onto the ideal of $\ZZ$-torsion elements in $H^*(BG;A)$. In particular, we have a short exact sequence 
$$
0\to {\rm Tor}_{\ZZ}(H^*(BG;A))\overset{j\beta^{-1}}{\to} \widehat{H}(\mathbf{B}_{\nabla}G;A)\overset{\mathcal{R}}{\to} I^*_A(\mathfrak{g}^{\vee})_{H}\to  0
$$
where $\beta:H^{\ast-1}(BG;A\otimes_{\ZZ}\CC^{\times})\to H^*(BG;A)$ is the Bockstein and ${\rm Tor}_{\ZZ}(H^*(BG;A))$ is the ideal of $\ZZ$-torsion elements. 
\end{theorem}
%In the case of $\mathbf{B}G$, we have the following. 
%\begin{theorem}
%The map 
% $j:H^{*}(\mathbf{B}G;A\otimes_{\ZZ}\CC^{\times})[-1]\to \widehat{H}^*(\mathbf{B}G;A)$
%is a surjective homomorphism of graded $A$-algebras, where $[-1]$ shifts the grading down by 1 and $H^{*}(\mathbf{B}G;A\otimes_{\ZZ}\CC^{\times})[-1]$ is equipped with the trivial algebra structure. If $G$ is compact and connected, then $j$ is an isomorphism and we have an identification 
%$$
%{\rm Tor}_{\ZZ}(H^*(BG;A))\cong \widehat{H}^*(\mathbf{B}G;A)
%$$
%where ${\rm Tor}_{\ZZ}(H^*(BG;A))$ is the graded submodule of $\ZZ$-torsion elements. 
%\end{theorem}

The second main result is a computation of the connective differential $K$-theory of the moduli stack $\mathbf{B}_{\nabla}G$. We have the following.

\begin{theorem}\label{theoremB}
Let $G$ be any Lie group and let $I_{\CC[u]}(\mathfrak{g}^{\vee})_k\subset I^*_{\CC[u]}(\mathfrak{g}^{\vee})$ denote the subring of $k$-integral elements. The curvature map defines a surjective homomorphism of graded rings
$$\mathcal{R}:\widehat{k}^*(\mathbf{B}_{\nabla}G)\to I^*_{\CC[u]}(\mathfrak{g}^{\vee})_{k}.$$
Moreover, the forgetful map $\mathcal{I}:\widehat{k}^*(\mathbf{B}_{\nabla}G)\to k^*(BG)$ maps the kernel of $\mathcal{R}$ isomorphically onto the image of the connecting homomorphism $\beta:k^{*-1}_{\CC^{\times}}(BG)\to k(BG)$. 
If $G$ is compact and connected, then $\mathcal{R}$ is isomorphism in even degrees $n\leq 0$. 
\end{theorem}

As a corollary, we show that the connective differential $K$-theory of $\mathbf{B}_{\nabla}G$ is isomorphic to just the topological $K$-theory of $BG$ in degree zero. Hence, by the Atiyah--Segal theorem, it is isomorphic to the completion of the representation ring at the augmentation ideal.

\begin{corollary}\label{corC}
Let $G$ be compact and connected. We have a commutative diagram of isomorphisms
$$
\xymatrix{
\widehat{k}^0(\mathbf{B}_{\nabla}G)\ar[r]^-{\mathcal{R}}\ar[d]^-{\mathcal{I}} & I_{\CC[u]}^0(\mathfrak{g}^{\wedge})_k
\\
R(G)^{\wedge}\ar[ru]_-{\rm ch} & 
}$$
where $R(G)^{\wedge}$ is the pro completion of $R(G)$ at the augmentation ideal. In particular, we have an isomorphism
$$\mathcal{I}:\widehat{k}^0(\mathbf{B}_{\nabla}G)\overset{\cong}{\to} R(G)^{\wedge}.$$

\end{corollary}
\def\smsset{{\rm C}^{\infty}{\rm s}\mathscr{S}{\rm et}}
\def\smset{{\rm C}^{\infty}\mathscr{S}{\rm et}}

\def\symgpd{{\rm s}\mathscr{G}{\rm pd}^{\otimes}}

\subsection*{Acknowledgments}

We thank Urs Schreiber, Arun Debray, and Severin Bunk for useful conversations about differential cohomology, which inspired the present work.

\def\smset{{\rm C}^{\infty}\mathscr{S}{\rm et}}
\section{Notation and conventions}

We will freely use the language of abstract homotopy theory. We will assume familiarity with model category theory, enriched model category theory, and Bousfield localization. Everything in the present work can almost certainly be done using Lurie's theory of $\infty$-categories \cite{HTT}, however it appears to us that  most calculations would involve several more steps and one would have to fix models for the quasi-categories anyways. For this reason, we have chosen to work in the 1-categorical setup from the start.

We will work with various presheaves on the site of cartesian spaces, whose definition is recalled below.
\begin{definition}
Let $(\Cart,J)$ denote the (1-categorical) Grothendieck site defined as follows.
\begin{itemize}
\item Objects are open subsets $U\subset \RR^n$ for some $n\in \NN$, such that $U$ is diffeomorphic to $\RR^n$. 
\item Morphisms are smooth maps between open subsets of $\RR^n$. 
\item The Grothendieck topology $J$ is generated by good open covers of cartesian spaces.
\end{itemize}
We call $(\Cart,J)$ the \emph{site of cartesian spaces}.
\end{definition}

\begin{notation}\label{notation}
We will use the following notation and conventions throughout.
\begin{enumerate}
\item Let ${\sf C}$ be a simplicial model category. Let $c,d\in {\sf C}$.  We denote the simplicial (nonderived) mapping space by  $\map(c,d)$. In general, if ${\sf C}$ is a ${\sf V}$-enriched model category, we denote the enriched (nonderived) hom by ${\sf V}(c,d)$, the only exception being the case of ${\sf V}=\mathscr{S}$.
\item
We denote the category of simplicial sets, equipped with the usual Kan--Quillen model structure, by $\mathscr{S}$. We call the objects of $\mathscr{S}$ \emph{spaces}.
\item We will denote the category of simplicial presheaves on $\Cart$, equipped with the local projective model structure, by $\sm\mathscr{S}$. Recall that the local model structure is obtained from the projective model structure by the left Bousfield localization at the following set of morphisms.
 \begin{itemize}
 \item Fix a cover $\{U_{a}\to U\}_{a\in A}$ of an object $U\in \Cart$. Let $\check{C}(\{U_{a}\})$ denote the simplicial presheaf in $\sm\mathscr{S}$ given by 
$$\check{C}(\{U_{a}\}):[n]\mapsto \coprod_{\alpha:[n]\to A} U_{\alpha},$$
where $U_{\alpha}=U_{\alpha_0}\cap U_{\alpha_1}\cap \hdots \cap U_{\alpha_n}$ (the $n$-fold intersections). The face and degeneracy maps are the obvious ones, given by inclusion and duplication of intersections. The localization is taken at the canonical morphisms of simplicial presheaves
\begin{equation}
\check{C}(\{U_{a}\})\to U,\label{cechmaps}
\end{equation}
induced by the the inclusions $U_{\alpha_0}\cap  U_{\alpha_1}\cap \hdots \cap U_{\alpha_n}\into U$. The resulting model category is simplicial, with 
$$\map(X,Y)_n=\set(\Delta^n,\map(X,Y))=\set(X\times \underline{\Delta}^n,Y),$$
where $\underline{\Delta}^n$ is the constant presheaf on $\Delta^n$. 

%More generally, if ${\sf V}$ is a closed symmetric monoidal category, then ${\sf C}$ can be enriched in ${\sf V}$ by ${\sf C}_{\sf V}(c,d):={\sf C}(c,d)\otimes {\bf 1}$, where ${\bf 1}$ is the monoidal unit. In this case, we let $\mathscr{S}{\rm h}_{(\infty,1)}({\sf C},{\sf V})$ denote the relative category with weak equivalences given by $S$-local weak equivalences, where $S$ is the set of morphisms of the form \eqref{cechmaps}, where now $U$ and $U_{\alpha}$ are embedded via the enriched Yoneda embedding. 
\end{itemize}
We call the objects of $\sm\mathscr{S}$ smooth stacks.

\item We will denote the category of simplicial symmetric spectra, equipped with positive stable model structure \cite{MMSS}, by $\Sp$.
%\item We will denote by $\Sp_{\geq 0}$ the category of connective spectra, equipped with local weak equivalences in the Bousfield--Friedlander model structure on $\Gamma$-spaces.
\item We will denote the category of presheaevs on $\Cart $ with values in $\Sp$, equipped with the local positive stable model structure, by $\sm\Sp$. The weak equivalences are given by $S$-local weak equivalences, where $S$ is the set of morphisms of the form
\begin{equation}\label{cechmorph}
\{F_n(\check{C}(\{U_{a}\}))\to F_n(U) \mid n\in \NN\}
\end{equation}
where $F_n$ is left adjoint to the evaluation functor ${\rm Ev}_n:\sm\Sp\to \sm\mathscr{S}$, which sends a presheaf of spectra $X$ to the presheaf of spaces given by taking its $n$-th level. Morally, the localization inverts the tensoring of {\v C}ech morphisms \eqref{cechmaps} with all desuspensions of the sphere (see \cite[Appendix B]{Dugger}).
\item We will denote the category of nonnegatively graded chain complexes by $\mathscr{C}{\rm h}_{\geq 0}$. The weak equivalences are quasi-isomorphisms. Similarly, the category of unbounded chain complexes with quasi-isomorphisms will be denotes $\mathscr{C}{\rm h}$. 
\item We will denote the category of presheaves of chain complexes on $\Cart$, equipped with the local projective model structure by $\sm\mathscr{C}{\rm h}_{\geq 0}$ and $\sm\mathscr{C}{\rm h}$, respectively. The localization is taken at the maps obtained by applying the normalized Moore functor $N$ (see the paragraph preceeding proposition \ref{DoldKan}) to the maps \eqref{cechmorph}.
\item We denote the canonical functors arising from the usual $t$-structure on $\mathscr{C}{\rm h}$ as follows. 
\begin{itemize}
\item The (stupid) truncation functors $\tau_{\geq 0},\tau_{\leq 0}:\mathscr{C}{\rm h}\to \mathscr{C}{\rm h}$  that sends a chain complex to the nonnegatively (nonpositively) graded chain complex 
\begin{align*}
\tau_{\leq 0}C &=(\cdots \to 0 \to C_0\to C_{-1}\to \cdots)
\\
\tau_{\geq 0}C &=(\cdots \to C_1 \to C_0\to 0\to \cdots)
\end{align*}
\item The (smart) truncation functor $z_{\geq 0}:\mathscr{C}{\rm h}\to \mathscr{C}{\rm h}$ that sends a chain complex to the nonnegatively graded chain complex 
$$z_{\geq 0}C=(\cdots \to C_1 \to Z(C_0)\to 0\to \cdots)$$
with cycles in degree $0$.
\item The (stupid) heart functor $T:\mathscr{C}{\rm h}\to \mathscr{C}{\rm h}$, $T=\tau_{\leq 0}\tau_{\geq0}$ and the (smart) heart functor $Z=\tau_{\leq 0}z_{\geq 0}$. 
\item For $n\in \ZZ$, the $n$-fold shift functor $s_n:\mathscr{C}{\rm h}\to \mathscr{C}{\rm h}$ that sends a chain complex to the shifted complex 
$$\hspace{1cm} s_n C=C[n] \quad C_{\bullet}[n]_i=C_{i-n}.$$
\end{itemize}
\end{enumerate}
\end{notation}

\section{Differential forms on smooth stacks}

In this section, we prove several results about differential forms on smooth stacks. These results are certainly known to experts, but we could not find a reference containing the level of generality needed for this work. We will begin with some basic definitions.

\begin{definition}
For each $q\in \NN$, we define the smooth abelian group of $q$-forms $\Omega^q$ on $\Cart$ as the presheaf 
$$U\mapsto \Omega^q(U)\qquad (f:U\to V)\mapsto f^*:\Omega^q(V)\to \Omega^q(U),$$
where $\Omega^q(U)$ denotes the group of complex-valued $q$-forms on $U$. 
Letting $q$-vary, we obtain a sheaf of unbounded nonpositively graded chain complexes
$$\Omega^*:=\left(\hdots \to 0\to \Omega^0\overset{d}{\to}\Omega^1\overset{d}{\to}\hdots \overset{d}{\to} \Omega^q\overset{d}{\to}\right),$$
where $d$ denotes the de Rham differential. Under the wedge product of forms, this is a sheaf of differentially graded commutative algebras. In particular, $\Omega^*$ is graded by $\ZZ_{\leq 0}$ by convention.
\end{definition}

More generally, we make the following definition. 
\begin{definition}
Let $A$ be a graded algebra of finite type over $\CC$ (i.e., finitely generated in each degree). We regard $A$ as a positively graded chain complex equipped with the trivial differential. The \emph{complex of forms with coefficients in $A$} is defined as the unbounded complex
$$
\Omega^*(-,A):=\Omega^*\otimes_{\CC}A, \quad d(\omega\otimes a)=d\omega\otimes a.
$$
The complex $\Omega^*(-,A)$ admits the structure of a presheaf of differentially graded algebras, with multiplication $\omega\otimes a \cdot \omega'\otimes a'=(-1)^{|\omega'||a|}\omega\wedge \omega'\otimes aa'$. 
\end{definition}

Since we include the case where $A$ is infinite dimensional, we might also want to consider the \emph{procomplete} tensor product $\Omega^*\widehat{\otimes }A$, defined as the limit of truncated complexes 
$$
 \Omega^*\widehat{\otimes} A=\lim_{i,j}\Omega^*_{\geq -i}\otimes A_{\leq j}.
$$ However, it turns out that the completed tensor product is naturally isomorphic to the usual tensor product in the case of differential forms. 
\begin{lemma}\label{procomplete}
Let $A$ be a differentially graded algebra of finite type over $\CC$.
The canonical map $\Omega^*\otimes A\to \Omega^*\widehat{\otimes}A$
is an isomorphism of presheaves of differentially graded algebras.
\end{lemma}
\proof
Fix $U\in \Cart$ and set $d={\rm dim}(U)$. Since the complex of forms $\Omega^*(U)$ vanishes in degrees below $-d$, the complex $\Omega^*(U)$ is bounded. The degree $n$ component of the canonical map is the map
$$\bigoplus_{j-i=n}\Omega^i(U)\otimes A_j\to \prod_{j-i=n}\Omega^i(U)\otimes A_j,$$
which on the $j,i$-summand sends $\omega\otimes a$ to the element whose $j,i$-component is $\omega\otimes a$ and all other components vanish. 
This is an isomorphism since there are only finitely many indices $i,j$ such that $j-i=n$ indexes a nonvanishing factor on the right. 
\endofproof

Next, we will discuss differential forms on smooth stacks. We begin with some preparation. Recall that the category of chain complexes admits a closed monoidal structure, with internal hom given by the complex
\begin{equation}\label{enrichedch}
\mathscr{C}{\rm h}(X,Y)_n:=\prod_{p\in \ZZ}\mathscr{A}{\rm b}(X_p,Y_{p+n}), \quad df=d_Yf-(-1)^nfd_X.
\end{equation}
The enrichment goes over to the category of smooth chain complexes $\sm\mathscr{C}{\rm h}$ (Notation \ref{notation}) over $\mathscr{C}{\rm h}$ in the usual way. Explicitly, the enriched hom is given by the same formula \eqref{enrichedch}, but with $\Ab(X_p,Y_{p+n})$ denoting the abelian group of natural transformations of presheaves of abelian groups (Notation \ref{notation}).

Given a smooth stack $X\in \sm\mathscr{S}$, we let $CX\in \sm\mathscr{C}{\rm h}$ denote the alternating face map complex of the corresponding free simplicial smooth abelian group. We make the following definition.
\begin{definition}\label{formsonstacks}
Let $X$ be a smooth stack and let $A$ be a graded algebra of finite type over $\CC$. We define the complex of differential forms on $X$ as the complex obtained via the enriched hom:
  $$
  \Omega^*(X,A):=\mathscr{C}{\rm h}(CX,\Omega^*(-,A)).
  $$
\end{definition}

Unwinding the definitions and using Lemma \ref{procomplete}, we compute the above complex in degree $n$ as
\begin{align*}
\Omega^n(X,A) &=\prod_{p\in \ZZ}\Ab(CX_p,\Omega^{*}(-,A)_{p-n}) 
\\
&\cong \prod_{p\in \ZZ}\Ab\left(CX_p,\prod_{j-i=p-n}\Omega^{i}\otimes A_j\right) 
\\
&\cong \prod_{p\in \ZZ}\prod_{j-i=p-n}\Omega^i(X_p)\otimes A_{j}
\\
&=\prod_p\prod_j\Omega^{j-p+n}(X_p)\otimes A_j,
\end{align*}
with differential given by 
$$
d(\omega\otimes a)=d\omega\otimes a-(-1)^n\delta \omega\otimes a, \quad \delta=\sum_{i=0}^n(-1)^id_i^*.
$$

%\begin{remark}
%The usual construction of the cup product using the Eilenberg--Zilber map gives the the complex $\Omega^*(X,A)$ a structure of a DGCA. More precisely, the product is given by the composition of the canonical map
%$$\mathscr{C}{\rm h}(CX,\Omega^*(-,A))\otimes \mathscr{C}{\rm h}(CX,\Omega^*(-,A))\to \mathscr{C}{\rm h}(CX\otimes CX,\Omega^*(-,A)\otimes \Omega^*(-,A))$$
%with the map 
%$$\mathscr{C}{\rm h}(CX\otimes CX,\Omega^*(-,A)\otimes \Omega^*(-,A))\to \mathscr{C}{\rm h}(CX,\Omega^*(-,A))$$
%that precomposes with the diagonal and Alexander--Whitney map $CX\to C(X\times X)\to CX\otimes CX$ and uses the product $\Omega^*(-,A)\otimes \Omega^*(-,A)\to \Omega^*(-,A)$. Explicitly, given elements 
%$$\omega\otimes a\cdot \omega'\otimes a'=\sum_{v+w=p+q}(-1)^{ij}d^*_{f,v}(\omega)\wedge d^*_{b,w}(\omega')\otimes aa' , $$
%where $d_{f,v}:X_{p+q}\to X_v$ is the front face map and $d_{b,w}:X_{p+q}\to X_w$ is the back face map.
%\end{remark}
\begin{remark}\label{procompstacks}
Note that if $A$ is finite dimensional, we have an isomorphism $\Omega^*(X,A)\cong \Omega^*(X)\otimes A$, which follows from the above computation by observing that there is a finite set of $j$'s indexing nonvanishing degrees and each $A_j$ is finitely generated. If $A$ is infinite dimensional, we instead obtain the procompleted tensor product. The two are not isomorphic for arbitrary stacks. 
\end{remark}
%\begin{remark}\label{cupprod}
%For any stack $X$, the complex of forms $\Omega^*(X,A)$ inherits the structure of a DGCA from $\Omega^*(-,A)$ in the usual way, via the composition of the maps
%\begin{align*}
%\iota:\mathscr{C}{\rm h}(CX,\Omega^*(-,A))\otimes \mathscr{C}{\rm h}(CX,\Omega^*(-,A)) &\to \mathscr{C}{\rm h}(CX\otimes CX,\Omega^*(-,A)\otimes \Omega^*(-,A))
%\\
%\Delta^*_{X,X}:\mathscr{C}{\rm h}(CX\otimes CX,\Omega^*(-,A)\otimes \Omega^*(-,A))&\to \mathscr{C}{\rm h}(C(X\times X),\Omega^*(-,A)\otimes \Omega^*(-,A))
%\\
%md^*:\mathscr{C}{\rm h}(C(X\times X),\Omega^*(-,A)\otimes \Omega^*(-,A))&\to \mathscr{C}{\rm h}(CX,\Omega^*(-,A))
%\end{align*}
%where $\iota(f\otimes f')(x\otimes y)=(-1)^{|f'||x|}f(x)\otimes f'(y)$, $\Delta^*_{X,X}$ is the Alexander--Whitney map, and $md^*$ is the composition of pullback by the diagonal and the multiplication $m:\Omega^*(-,A)\otimes \Omega^*(-,A)\to \Omega^*(-,A)$. 
%\end{remark}
 A priori, it may not be clear that Definition \ref{formsonstacks} is the correct one, since the homs are not derived. However, in good cases it can be shown to be derived. For example, if one takes $X=M$ to be a smooth manifold, then the derived complex can be computed by the strict one, which is essentially the content of the {\v C}ech--de Rham theorem for smooth manifolds.

The next proposition generalizes the {\v C}ech--de Rham theorem to a larger class of smooth stacks. We will begin with some preparation. Let $\mathcal{A}\subset \sm\mathscr{S}{\rm et}$ denote the class of objects that are jointly acyclic for the functors $\Omega^q$ with $q>0$. That is, for a fixed cofibrant replacement $QS$ of $S\in \smset\subset \sm\mathscr{S}$ in the local projective model structure, an object $S\in \mathcal{A}$ if and only if the cohomology of the complex
\begin{equation}\label{acyclic}
\Omega^q(QS_p), \qquad d=\sum_{i=0}^p(-1)^id_i^*
\end{equation}
vanishes for all $p>0$ and $q\geq 0$. This is equivalent to the vanishing of the higher Ext groups ${\rm Ext}^p_{\sm\mathscr{A}{\rm b}}(F(S),\Omega^q)$ for $p>0$ and $q\geq 0$, where $F:\smset\to \sm\mathscr{A}{\rm b}$ is the free functor. This is further equivalent, using the simplicially enriched Dold-Kan correspondence, to the vanishing of the homotopy groups of the derived mapping space $\pi_pR\map(S,\Omega^q)$ for $p>0$ and $q\geq 0$. 

 \begin{proposition}
Let $A$ be a graded $\CC$-algebra of finite type. Let $Q$ be a cofibrant replacement functor in $\sm\mathscr{S}$. Let $X$ be a simplicial object in $\mathcal{A}$. 
Then the map $QX\to X$ induces a quasi-isomorphism of complexes:
  $$
  \Omega^*(X,A)\overset{\simeq}{\to} \Omega^*(QX,A).
  $$
  \end{proposition}
  \proof
 The complexes $\Omega^*(X,A)$ and $\Omega^*(QX,A)$ are given by the enriched homs
\begin{align*}
& \Omega^*(X,A)=\mathscr{C}{\rm h}(CX,\Omega^*(-,A))
\\
& \Omega^*(QX,A)=\mathscr{C}{\rm h}(CQX,\Omega^*(-;A)).
\end{align*}
We claim that the map $QX\to X$ induces a quasi-isomorphism 
$$\mathscr{C}{\rm h}(CX,\Omega^*(-,A))\to \mathscr{C}(CQX,\Omega^*(-,A)),$$
which will prove the claim. First, observe that $(QX)_p$ is acyclic for $\Omega^q$. Indeed, by the characterization of cofibrant  objects in $\sm\mathscr{S}$ \cite[Corollary 9.4]{Dug}, $(QX)_p$ is a disjoint union of representables, hence cofibrant when regarded as an object in $\sm\mathscr{S}$. Since $\Omega^q$ forms a sheaf of abelian groups, it is fibrant in the local model structure on $\sm\mathscr{S}$. Therefore, the mapping space $\map(QX_p,\Omega^q)\cong N\set(QX_p,\Omega^q)$ is already derived. 
%$$
%H^i(C(Q(X))_p;\Omega^q)\cong \prod_{U_p\to \hdots U_0\to X_p}H^i(U;\Omega^q)\cong 0, \quad\forall i>0
%$$

The class of sheaves of abelian groups that are acyclic for $\Omega^q$ is closed under direct sums. Hence, the mapping cone of the quasi-isomorphism $CQX\to CX$ is levelwise acyclic for $\Omega^q$ and is quasi-isomorphic to $0$. Since $\Omega^q$ preserves exact complexes of acyclics, it follows by the long exact sequence that we have a quasi-isomorphism
 $$
 \mathscr{C}{\rm h}(CX,\Omega^q)\to \mathscr{C}{\rm h}(CQX,\Omega^q), \quad \forall q\geq 0.
 $$
Finally, since $A$ is of finite type, we can write $\Omega^*(-,A)$ levelwise as a direct sum of $\Omega^q$, which by the above implies that we have a quasi-isomorphism
 $$\mathscr{C}{\rm h}(CX,\Omega^*(-,A))\to \mathscr{C}{\rm h}(CQX,\Omega^*(-,A)).$$
\endofproof

The {\v C}ech cohomology of a smooth stack $X\in \sm\mathscr{S}$ with coefficients in a graded ring $A$ can be defined via the derived hom
$$H_{\check{C}}^n(X;A)=\pi_0\map(QX,K(\underline{A}[n])),$$
where $QX$ is a cofibrant replacement in the local model structure, and $K(\underline{A}[n])$ is the locally constant sheaf obtained by applying the Dold--Kan functor in Proposition \ref{DoldKan} to the sheaf of chain complexes given by locally constant $A$-valued functions $\underline{A}$, shifted in degree, equipped with the trivial differential. Equivalently, using the enriched Dold--Kan correspondence, the above is just the cohomology of the complex $\mathscr{C}{\rm h}(CQX,\underline{A})$. If $X=M$ is a smooth manifold and the cofibrant replacement is given by the {\v C}ech resolution by a good open cover, then the above agrees with the usual definition of {\v C}ech cohomology of a smooth manifold, with coefficients in $A$.

%\begin{remark}\label{cupcech}
%The {\v C}ech cohomology of a smooth stack $X$ with coefficients in $A$ admits the structure of a graded algebra over $A$ via the cup product. The product is defined the same was as in the case of forms (Remark \ref{cupprod}), replacing $X$ with $QX$ and $\Omega^*(-,A)$ with $\underline{A}$.
%\end{remark}

\begin{definition}
Let $X\in \sm\mathscr{S}$ be a smooth stack and let $QX$ be a cofibrant replacement. Let $A$ be a graded $\CC$-algebra of finite type. We define the de Rham cohomology of $X$ with coefficients in $A$ as
$$H^*_{\rm dR}(X;A):=H_{-\ast}(\Omega^*(X,A)).$$
We define the {\v C}ech cohomology with coefficients in the sheaf of locally constant $A$-valued functions $\underline{A}$ by 
$$H^*_{\check{C}}(X;\underline{A}):=H_{-\ast}( \mathscr{C}{\rm h}(CQX,\underline{A})).$$
%Both admit the structure of graded commutative algebras over $A$ via the products in Remark \ref{cupcech} and \ref{cupprod}.
\end{definition}

%\begin{definition}
%Let $A$ be a graded commutative algebra over $\CC$, Let $X$ be a smooth stack. We define the de Rham cohomology of $X$ with coefficients in $A$ in degree $n\geq 0$ as
%$$H^n_{{\rm dR}}(X;A):=H_n(\Omega^*(X,A[n])).$$
%\end{definition}
%
%\begin{remark}
%Both $H^*_{\check{C}}(X;A)$ and $H_{dR}^*(X;A)$ admit the structure of graded commutative algebras over the graded algebra $A$ \cite{??} in the obvious way. Explicitly, given homogeneous elements $x=\sum_{i} x_i, y=\sum_jy_j, x_i\in H^{n+i}(X;A_i), y_j\in H^{m+j}(X;A_j)$, the product is given by 
%$$x\cup y=\sum_{i,j}x_i\cup y_j, \quad x_i\cup y_j\in H^{n+m+i+j}(X;A_{i+j})\subset H^{n+m}(X;A).$$
%
%\end{remark}

The following corollary is a generalized version of the {\v C}ech--de Rham theorem for levelwise acyclic smooth stacks.

\begin{corollary}
Let $X$ be a simplicial object in $\mathcal{A}$ and let $A$ be a graded $\CC$-algebra of finite type. Then we have an isomorphism of graded rings
$$H^*_{\rm dR}(X;A)\cong H_{\check{C}}^*(X;A).$$
\end{corollary}
\proof
By the Poincar\'e lemma, the inclusion $\underline{A}\into \Omega^*(-,A)$ is a quasi-isomorphism of sheaves of chain complexes. Since $CQX$ is cofibrant in the projective model structure on $\sm\mathscr{C}{\rm h}$, the induced map 
$$\mathscr{C}{\rm h}(CQX,\underline{A})\to \mathscr{C}{\rm h}(CQX,\Omega^*(-,A))=\Omega^*(QX,A)$$
is also a quasi-isomorphism. By the previous proposition, we have a quasi-isomorphism 
$\Omega^*(X,A)\to \Omega^*(QX,A)$, which proves that we have an isomorphism on cohomology. Since the inclusion $\underline{A}\into \Omega^*(-,A)$ is a map of presheaves of differentially graded algebras, it follows that the induced map on cohomology is a map of graded rings.
\endofproof

\begin{corollary}
The {\v C}ech--de Rham theorem holds for simplicial manifolds and homotopy invariant stacks.
\end{corollary}

We now specialize to the case where $X$ is the moduli stack of smooth principal $G$-bundles or principal $G$-bundles with connection. Recall that, on the site of cartesian spaces, the moduli stack of principal $G$-bundles $\mathbf{B}G$ is given by taking the nerve of the action groupoid corresponding to the trivial action of the smooth group $C^{\infty}(-,G)$ on the terminal object $\ast$. The moduli stack of principal $G$-bundles with connection $\mathbf{B}_{\nabla}G$ is given by taking the nerve of the action groupoid corresponding to the action of the smooth group $C^{\infty}(-,G)$ on the smooth set of $\mathfrak{g}$-valued 1-forms $\Omega^1(-,\mathfrak{g})$ by gauge transformations.

\begin{notation}\label{oppgr}
In order to get the gradings to work out correctly, we will use the following notation throughout. Let $A$ be a $\ZZ$-graded algebra. In particular, any we can regard any $\NN$-graded algebra as a $\ZZ$-graded algebra that is trivial in negative degrees. The sign involution on $\ZZ$ gives rise to an involution on $\ZZ$-graded algebras, which we denote by 
$${\rm gr}_{\ZZ}\mathscr{A}{\rm lg}\to {\rm gr}_{\ZZ}\mathscr{A}{\rm lg}, \quad A\mapsto \overline{A}.$$
Explicitly, we have $\overline{A}_{-n}=A_n$. In particular, a nonnegatively graded algebra is sent to a nonpositively graded algebra and vice versa. 
\end{notation}

\begin{definition}\label{invcoA}
Let $A$ be a graded $\CC$-algebra of finite type. Let $I^*(\mathfrak{g}^{\vee})={\rm Sym}^*(\mathfrak{g}^{\vee})^G$ denote the symmetric algebra of Ad-invariant polynomials on $\mathfrak{g}^{\vee}$, graded by twice the degree. Using Notation \ref{oppgr}, let
$$I_{A}^*(\mathfrak{g}^{\vee}):=I^*(\mathfrak{g}^{\vee})\widehat \otimes_{\CC} \overline{A}$$ denote the pro-completed tensor product of corresponding $\ZZ$-graded algebras (i.e., the pro-completion is taken levelwise). We call $I^*_{A}(\mathfrak{g}^{\vee})$ the \emph{$\ZZ$-graded algebra of invariant polynomials with coefficients in $A$}. If $A$ is an graded ring of finite type, we define $I_{A}^*(\mathfrak{g}^{\vee}):=I^*_{A\otimes_{\ZZ} \CC}(\mathfrak{g}^{\vee})$.
\end{definition}

As an immediate corollary of the {\v C}ech--de Rham theorem, we obtain the Chern--Weil theorem with coefficients in an arbitrary graded algebra of finite type.
\begin{corollary}\label{chernweil}
Let $G$ be a compact and connected Lie group. Let $A$ be a graded $\CC$-algebra of finite type. We have isomorphisms of graded rings
\begin{equation}\label{bottcderham}
H^*(BG;A)\cong H^*_{\check{C}}(\mathbf{B}_{\nabla}G,A)\cong H_{\check{C}}^*(\mathbf{B}G;A)\cong I^{*}_A(\mathfrak{g}^{\vee}).
\end{equation}
\end{corollary}
\proof
Since $\mathbf{B}G$ is a simplicial manifold, the {\v C}ech--de Rham theorem holds. By the main theorem of \cite{Bott}, if $G$ is compact then the cohomology of the complex of differential forms $\Omega^*(\mathbf{B}G)$ is isomorphic to the graded algebra of Ad-invariant polynomials $I^*(\mathfrak{g}^{\vee})$, graded by twice the degree. According to our grading convention, this means
$$H_{-\ast}(\Omega^*(\mathbf{B}G,A))\cong I^*(\mathfrak{g}^{\vee})$$
(see Notation \ref{oppgr}). By Lemma \ref{procomplete} (and remark \ref{procompstacks}) we have isomorphisms
\begin{align*}
\Omega^*(\mathbf{B}G,A) &\cong \mathscr{C}{\rm h}(C(\mathbf{B}G),\Omega^*\otimes A)
\\
&\cong \mathscr{C}{\rm h}(C(\mathbf{B}G),\Omega^*\widehat{\otimes}A)
\\
&\cong \varprojlim \mathscr{C}{\rm h}(C(\mathbf{B}G),\Omega^*_{\geq -j}\otimes A_{\leq k})
\\
&\cong \varprojlim \Omega^*(\mathbf{B}G)_{\geq -j}\otimes A_{\leq k}.
\end{align*}
Taking homology and using the fact that $A$ is equipped with the trivial differential, we compute
$$H_{-\ast}(\Omega^*(\mathbf{B}G,A))\cong H_{-*}\left(\varprojlim \Omega^*(\mathbf{B}G)_{\geq -j}\otimes A_{\leq k}\right)\cong (I^*(\mathfrak{g}^{\vee})\widehat{\otimes} \overline{A})_{*}=I^{*}_{A}(\mathfrak{g}^{\vee}),$$
where we have used the fact that for all $q\in \ZZ$, the tower of groups 
$$H_{-q}(\Omega^*(\mathbf{B}G)_{\geq - j}\otimes A_{\leq k})\cong (I^*(\mathfrak{g}^{\vee})_{\leq j}\widehat{\otimes} \overline{A}_{\geq -k})_q$$ satisfies the Mittag-Leffler condition, so that $\lim^1$ vanishes.  

By the Poincar\'e  lemma, we have a quasi-isomorphism of complexes $\underline{A}\into \Omega^*(-,A)$. Then the {\v C}ech--de Rham theorem implies the right isomorphism in \eqref{bottcderham}. The remaining isomorphisms are obtained as follows. Let $\underline{A}[n]$ denote the $n$-fold shift of $\underline{A}$ and let $K(\underline{A}[n])$ denote the smooth simplicial set obtained via the Dold--Kan correspondence (Proposition \ref{DoldKan}). Observe that $K(\underline{A}[n])$ is fibrant in the  projective model structure, since it is an objectwise Kan complex. It is also satisfies homotopy descent on cartesian spaces (see, for example, the proof of \cite[Proposition 3.3.9]{Urs}). Hence it is fibrant in the local projective model structure.  It follows that the derived mapping space can be computed by $R\map(Q\mathbf{B}G,K(\underline{A}[n]))$, whose connected component give the {\v C}ech cohomology in degree $n$.  

On the other hand, the derived mapping space can be computed using the Quillen adjunction $|\cdot|\dashv \delta$ (Proposition \ref{cohesive}), which implies that the derived mapping space $R\map(\mathbf{B}G,K(\underline{A}))$ is weakly equivalent to $R\map(|\mathbf{B}_{\nabla}G|,K(A))$. Since $|\mathbf{B}_{\nabla}G|$ is weakly equivalent to $BG$ (Example \ref{clssbg}), the left isomorphism is induced by the corresponding weak equivalence on derived mapping spaces. Similarly, the middle isomorphism follows from the weak equivalence $|\mathbf{B}G\vert\simeq BG$ (Example \ref{clssbg}). 

Finally, since the ring structure is induced from the coefficient ring $A$, and the inclusion $\underline{A}\into \Omega^*(-,A)$ is a map of graded algebras, it follows that all the isomorphisms \eqref{bottcderham} are isomorphisms of rings. 
\endofproof

The following example will be useful in the computation of the differential $k$-theory of $\mathbf{B}_{\nabla}G$.
\begin{example}\label{invariantpolys}
Let $G$ be compact and connected, with Lie algebra $\mathfrak{g}$. Fix a maximal torus $T\subset G$ with Cartan algebra $\mathfrak{t}$ and let $\{\lambda_i\}_{i=1}^n$ be a basis of $\mathfrak{t}^{\vee}$.  Let $W(G)$ denote the Weyl group. By the Chevalley restriction theorem, we have an isomorphism of graded polynomial algebras $I^*(\mathfrak{g}^{\vee})\cong \CC[\lambda_1, \hdots, \lambda_n]^{W(G)}$, $|\lambda_i|=2$. Let $A=\CC[u]$, $|u|=2$. Then we compute the pro-completed tensor product as 
$$I^*_{\CC[u]}(\mathfrak{g}^{\vee})=\varprojlim I^*(\mathfrak{g}^{\vee})_{\leq j}\otimes \overline{\CC[u]}_{\leq k}=\varprojlim \CC[\lambda_1, \hdots, \lambda_n]^{W(G)}/J^j\otimes \CC[\bar u]/(\bar u^k),$$
where $J$ is the augmentation ideal and $|\bar u|=-2$. By invariant theory, the augmentation ideal $J$ admits a minimal set of $n$ algebraically independent generators $I_1,\hdots,I_n$, and we have an identification of $W(G)$-invariants as the polynomial algebra
$$\CC[\lambda_1,\hdots,\lambda_n]^{W(G)}\cong \CC[I_1,\hdots,I_n].$$
To describe the ring structure of $I^*_{\CC[u]}(\mathfrak{g}^{\vee})$, it is useful to consider the canonical map of graded algebras $\CC[\bar u]\to \CC[\bar u,\bar u^{-1}]$, $\bar u\mapsto \bar u$. This map produces an injective map
$$I^*(\mathfrak{g}^{\vee})\widehat{\otimes}\CC[\bar u]\to I^*(\mathfrak{g}^{\vee})\widehat{\otimes} \CC[\bar u,\bar u^{-1}]\cong I(\mathfrak{g}^{\vee})^{\wedge}[v,v^{-1}],$$
where $I(\mathfrak{g}^{\vee})^{\wedge}$ is the completion (in the category of nongraded rings) at $J$. The isomorphism on the right sends $1\mapsto 1$, $\bar u\mapsto v^{-1}$ and $I_j\mapsto v^{|I_j|}I_j$. The composite map is injective, with image generated by $1,v^{-1},v^{|I_j|}I_j$. This subring appears to be related to the (completed) modified Rees ring of the representation ring $R(G)$,  which was used in \cite{Greenlees} to study the coefficient ring of equivariant connective $K$-theory $k^*_{G}$, although we will not pursue the precise relationships here.
\end{example}

%
%\begin{remark}
%The smooth set $\Omega^1(-;\mathfrak{g})$ is not acyclic for $\Omega^q$. One can find a nontrivial {\v C}ech cocycle in degree 1 as follows. The cofibrant replacement of $\Omega^1(-;\mathfrak{g})$ is the simplicial presheaf that in degree $n$ is given by 
%$$\coprod_{U_{n}\to\hdots \to U_0\to \Omega^1(-;\mathfrak{g})}U_n$$
%The face and degeneracy maps are the obvious ones that remove or duplicate the $i$-th component in the chain. Mapping into the presheaf of simplicial abelian groups $\Omega^q[1]$ and using the Dold-Kan correspondence, one can compute $H^1(\Omega(-;\mathfrak{g}),\Omega^q)$ as $H^1$ of the cochain complex 
%$$\prod_{U_n\to \hdots U_0\to \Omega^1(-;\mathfrak{g})}\Omega^q(U_n), \quad d=\sum_j(-1)^jd_j^*$$
%Thus, a closed element in degree 1 is a collection of forms 
%$$\omega_{f:U_1\to U_0,\mathcal{A}\in \Omega^1(U_0;\mathfrak{g})}\in \Omega^q(U_1)$$ 
%such that for every composable pair $fg:U_2\to U_1\to U_0$, we have 
%$$g^*\omega_{f,\mathcal{A}}=\omega_{fg,\mathcal{A}}-\omega_{g,f^*\mathcal{A}}$$
%Let $p_q\in \mathcal{I}(G)^q$ be an invariant polynomial of degree $q$. Define 
%$$\omega_{f,\mathcal{A}}={\rm CS}(\mathcal{A}).$$
%Then 
%$$g^*\omega_{f,\mathcal{A}}=g^*e^{f^*\mathcal{A}}$$
%and 
%$$\omega_{fg,\mathcal{A}}-\omega_{g,f^*\mathcal{A}}=j_{fg}^*\mathcal{A}_{fg}-j_g^*f^*\mathcal{A}$$
%Hence, $\omega_{f,\mathcal{A}}$ is closed. If it were exact, then there would be $\tau_{\mathcal{A}}$ such that 
%$$\omega_{f,\mathcal{A}}=\tau_{f^*\mathcal{A}}-f^*\tau_{\mathcal{A}}$$
%
%\end{remark}

\begin{definition}
Let $X$ be a smooth stack. Let $A$ be a graded $\CC$-algebra of finite type. We define the complex of basic forms on $X$ as the equalizer 
$$\Omega^*(X,A)_{bas}:=\lim\left\{\xymatrix{
\Omega^*(X_0,A)\ar@<.1cm>[r]^-{d_0^*}\ar@<-.1cm>[r]_-{d_1^*}  & \Omega^*(X_1,A)}\right\}$$
We say that that a basic form $\omega$ is closed if $d\omega=0$, where $d$ is the differential inherited from $\Omega^*(-,A)$. We denote the group of basic closed forms by $\Omega^*(X,A)_{bas,cl}$.
\end{definition}

When $M$ is a smooth manifold equipped with a $G$-action, applying the previous definition to the action groupoid $M\rtimes G$ agrees with the usual notion of basic form from equivariant de Rham theory. Indeed, any form in the equalizer is manifestly $G$-invariant. For any $\xi\in \mathfrak{g}$, any form in the equalizer is also annihilated by contraction with the induced vector field $\xi^{\sharp}\in \Gamma(M,TM)$, which follows by observing that the induced vector field $\xi^{\sharp}$ on $G\times M$ (with $G$ acting on the first factor) is given by  $\xi^{\sharp}(g,x)=(g\xi,0)$ and this is mapped to $\xi^{\sharp}$ by the differential of the action. Hence, the differential of $p$ maps this vector field to zero and $a^*i_{\xi^{\sharp}}\omega=i_{\xi^{\sharp}}a^*\omega=i_{\xi^{\sharp}}p^*\omega=0$.

The following proposition will be very useful in computing differential forms on smooth stacks. 

\begin{proposition}
Let $X$ be a smooth stack. Let $A$ be a graded $\CC$-algebra of finite type. Let $n$ be a nonnegative integer. Let $Z$ be the smart heart functor (Notation (ix) \ref{notation}) and let $[n]$ denote the $n$-fold shift. Let $Q$ be a cofibrant replacement functor on $\sm\mathscr{S}$. We have an isomorphism 
$$\Omega^n(X,A)_{bas,cl}\overset{\simeq}{\to}Z(\Omega^*(QX,A[n]))=\Omega^n_{cl}(QX,A).$$
\end{proposition}
\proof
The smooth set  $Z(\Omega^*(-,A[n]))=\Omega_{cl}^n(-,A)$ is an objectwise Kan complexes (being zero truncated). It is also local with respect to {\v C}ech covers, since closed differential forms form a sheaf. Hence, $\Omega_{cl}^n(-,A)$ is a fibrant object in the local projective model structure. 
Since $\Omega^n_{cl}(-,A)$ is zero truncated for all $n\in \NN$, Lemma \ref{truncate} applies and the canonical map $QX\to \pi_0(QX)$ induces a weak equivalence
$$N\set(\tilde \pi_0(X),\Omega_{\rm cl}^n(-,A))\to \map(QX,Z(\Omega^*(-,A[n]))),$$
where $N:\set\to\mathscr{S}$ is the nerve.
But the left side is manifestly in natural bijection with the nerve of the set of basic forms. This completes the proof.
\endofproof

%\begin{remark}
%There is an obvious comparison map to the complex defined by Hopkins and Freed, which is given as follows. The truncation functor $\tau_0:\sset\to \set$ is left adjoint to constant diagram functor. The adjunction unit evaluated on (the nerve of) the action groupoid $S\rtimes G$ gives a map 
%$$S\times G\to \tau_0(S\rtimes G)$$
%This induces an induced map on double complexes 
%$$ \Omega^q(S)_b=\Omega^q(\tau_0(S\rtimes G))\to \Omega^q(S\times G^p),$$
%hence an induced map on totalizations. This yields an inclusion of $\Omega^*(S)_b$ as a subcomplex of the total complex of equivariant forms on $S$, which in degree $n=p+q$ corresponds to the summand indexed by $(0,q)$.
%\end{remark}

\begin{corollary}\label{basicfrmsbg}
 Let $A$ be a graded $\CC$-algebra of finite type. Let $n$ be a nonnegative integer. Let $Z$ be the smart heart functor (Notation (ix) \ref{notation}) and let $[n]$ denote the $n$-fold shift. Let $Q$ be a cofibrant replacement functor on $\sm\mathscr{S}$. We have an isomorphism
\begin{equation}\label{fhbasic}
I^n_{A}(\mathfrak{g}^{\vee})\cong \Omega^n(\mathbf{B}_{\nabla}G,A)_{bas,cl}\to Z(\Omega^*(Q(\mathbf{B}_{\nabla}G),A[n]))
\end{equation}
\end{corollary}
\proof
Setting $X=\mathbf{B}_{\nabla}G$ in the previous corollary gives the isomorphism appearing as the right arrow in \eqref{fhbasic}. The isomorphism on the left of \eqref{fhbasic} follows from the work of Freed--Hopkins as follows. The main theorem of \cite{FH} gives an isomorphism 
$$\set(\tilde \pi_0(\mathbf{B}_{\nabla}G),\Omega^n)\cong \Omega^n(\mathbf{B}_{\nabla}G)_{bas,cl}\cong I^{n}(\mathfrak{g}^{\vee}),$$
%where the overline (Notation \ref{oppgr}) comes from our choice of grading on forms. In particular,
%$$\Omega^n(\mathbf{B}_{\nabla}G)=(\Omega^*(\mathbf{B}_{\nabla}G)_{bas,cl})_{-n}=I^n(\mathfrak{g}^{\vee}).$$ 
To get the result for arbitrary coefficients, observe that the calculation after Definition \ref{formsonstacks} computes the degree $n$ component of the procompleted tensor product (Remark \ref{procompstacks}). Setting $X=\tilde \pi_0 \mathbf{B}_{\nabla}G$ there yields 
\begin{align*}
\Omega^n(\mathbf{B}_{\nabla}G,A)_{bas,cl} &=\Omega^n_{cl}(\tilde \pi_0(\mathbf{B}_{\nabla}G),A)
\\
& =\prod_{j}I^{j+n}(\mathfrak{g}^{\vee})\otimes \overline{A}_{-j}=I^{n}_A(\mathfrak{g}^{\vee}).
\end{align*}

%The composite functor $KZs_n:\sm\mathscr{C}{\rm h}\to \smset$ (See (ix) in Notation \ref{notation}) sends a sheaf of chain complexes to the smooth set obtained by shifting up in degree by $n$, throwing away all degrees but zero, and puts cycles in degree zero. This functor is right adjoint, with corresponding left adjoint given by $s_{-n}iF:\smset\to \sm\mathscr{C}{\rm h}$, where $F$ is the free functor and $i:\sm\mathscr{A}{\rm b}\to \sm\mathscr{C}{\rm h}$ is the inclusion, and $s_{-n}$ is the downward shift functor. Hence, $KZs_n$ preserves limits. We have
%\begin{align*}
%\Omega^n(\mathbf{B}_{\nabla}G,A)_{bas,cl} &\cong \sm\set(\tilde \pi_0(\mathbf{B}_{\nabla}G),KZ(\Omega^*\otimes A[n]))
%\\
%&\cong \sm\set(\tilde \pi_0(\mathbf{B}_{\nabla}G),KZ(\Omega^*\hat{\otimes}A[n]))
%\\
%&\cong \varprojlim \sm\set(\tilde \pi_0(\mathbf{B}_{\nabla}G),KZ(\Omega^*_{\leq j}\otimes A_{\leq k}[n]))
%\\
%&\cong \varprojlim (\Omega^*_{\leq j}(\mathbf{B}_{\nabla}G)_{bas,cl}\otimes A_{\leq k})_n
%\\
%&\cong \varprojlim (I^*(\mathfrak{g})_{\leq j}\otimes A_{\leq k})_n=: I_A^*(\mathfrak{g}^{\vee})_n.
%\end{align*}
%where going from the 3rd to the 4th line is a simple unwinding of the definitions, using the identification in \ref{?}.
\endofproof

\section{The differential cohomology of $\mathbf{B}_{\nabla}G$}

We are ready to prove the first of the main theorems. We first recall some basic definitions. Let $A$ be a graded ring of finite type. Assume that $A$ is levelwise torsion free (as a $\ZZ$-module). The exponential sequence $\ZZ \into \CC\to \CC^{\times}$ induces an exact sequence of graded $A$-modules
$$0\to A\to A\otimes_{\ZZ}\CC\to A\otimes \CC^{\times}\to 0$$
and hence a long fiber/cofiber sequence of corresponding Eilenberg--MacLane spectra 
\begin{equation}\label{bockstien}
\hdots \to \Sigma^{-1}H(A\otimes_{\ZZ}\CC^{\times})\overset{\beta}{\to} HA\to H(A\otimes_{\ZZ}\CC)\to H(A\otimes_{\ZZ}\CC^{\times})\overset{\beta}{\to} \Sigma HA\to \hdots 
\end{equation}
The map $\beta$ is the Bockstein map. It induces a natural map on cohomology with coeffcients in $A$: $\beta:H^*(-;A\otimes_{\ZZ}\CC^{\times})\to H^{*+1}(-;A)$. 
%All maps are compatible with the $A$-module structure, since the original exact sequence is exact in the category of graded $A$-modules. 

\begin{remark}
As a model for the Eilenberg--MacLane spectrum functor, we use the \emph{naive} Eilenberg--MacLane spectrum functor that lands in sequential spectra. This is necessary to compare with symmetric module spectra, through a zig-zag of Quillen equivalences (\cite[Appendix B]{ShipSchwede}). In this section, we take $\Sp$ to be the category of sequential spectra with the stable model structure. 
\end{remark}
Following \cite{BNV}, we have the following definition.

\begin{definition}
Let $A$ be a graded ring of finite type. We define the differential cohomology spectrum with coefficients in $A$ as the homotopy fiber product in $\sm\Sp$:
\begin{equation}\label{diffcoh}
\xymatrix{
\widehat{H}A(n)\ar[r]^-{\mathcal{R}}\ar[d]^-{\mathcal{I}} & H(Z\Omega^*(-,A[n]))\ar[d]
\\
\Sigma^nH\underline{A}\ar[r] & H(\Omega^*(-,A[n])),
}
\end{equation}
where the map bottom horizontal map is induced by the canonical inclusion $\underline{A}[n]\to \Omega^*(-,A[n])$ as locally constant functions. The right vertical map is induced by the the canonical inclusion $Z\Omega^*(-,A[n])\to \Omega^*(-,A[n])$ of cycles in degree zero.

The maps $\mathcal{R}$ and $\mathcal{I}$ induce maps on cohomology, which we call the \emph{curvature map} and \emph{forgetful map}, respectively.
\end{definition}

The above definition of differential cohomology differs slightly from the one in \cite{BNV}. In particular, we have chosen to work with the zero truncated complex $Z\Omega^*(-,A[n])$, which allows us to apply Corollary \ref{basicfrmsbg} to compute forms on $\mathbf{B}_{\nabla}G$. Nevertheless, the sheaf of spectra \eqref{diffcoh} is equivalent to the one in \cite{BNV}. Indeed, the inclusion $Z\Omega^*(-,A[n])\into \tau_{\leq 0}\Omega^*(-,A[n])$ is a quasi-isomorphism of shevaes complexes, by the Poincar\'e lemma. Similarly, the inclusion $\underline{A}[n]\otimes \CC\into \Omega^*(-,A[n])$ is also quasi-isomorphism. Hence the homotopy fiber product in \eqref{diffcoh} is weakly equivalent to the one considered in \cite{BNV}. 

\begin{lemma}\label{bnlemma}
Let $\mathbf{B}^n_{\nabla}A$ denote the $0$-th layer of the spectrum $\widehat{H}A(n)$. Then $\mathbf{B}^n_{\nabla}A$ fits into a homotopy cartesian square in the $\sm\mathscr{S}$:
\begin{equation}\label{hsquarepull}
\xymatrix{
\mathbf{B}^n_{\nabla}A\ar[r]^-{\mathcal{R}}\ar[d]^-{\mathcal{I}} & \Omega^n_{\rm cl}(-,A) \ar[d]
\\
K(\underline{A}[n])\ar[r] & K(\Omega^*(-,A[n]))
}
\end{equation}
\end{lemma}
\proof
Since Eilenberg--MacLane spectra are $\Omega$-spectra and these are fibrant in the stable model structure $\Sp$, it follows that the smooth spectra appearing in the three corners of \eqref{diffcoh} are fibrant in the projective model structure on $\sm\Sp$. We claim that they are also fibrant in the local model structure. Indeed, by definition of the localization in Notation \ref{notation}, local objects are precisely those spectra which satisfy descent levelwise. By definition, the $j$-th level of the Eilenberg--MacLane spectrum is given by applying $K$ to the the $j$-fold shift of the corresponding complex. For $j>0$, the relevant complexes are levelwise direct sums of differential forms and locally constant sheaves. These complexes all satisfy homotopy descent, which can be seen be computing the relevant homotopy limits as totalizations and using the fact the forms are acyclic and that cartesian spaces are contractible. We leave the details to the reader. For $j=0$, we have  $H(Z\Omega^*(-,A[n]))_0=\Omega^n_{cl}(-,A)$. Since this presheaf is zero truncated, homotopy descent reduces to just the ordinary sheaf property, which is satisfied for closed forms. 

 With these observations, we see that to compute the homotopy pullback in $\sm\Sp$, it follows from \cite[Proposition 3.3.16]{Hirschhorn} that it suffices to replace one of the legs by a fibrant resolution in just the projective model structure, hence a levelwise fibration, and compute the levelwise strict pullback. Computing the homotopy pullback this way shows in particular that the zeroth level of $\widehat{H}A(n)$ fit into homotopy cartesian square \eqref{hsquarepull}, as claimed.
\endofproof
%Since $A$ is a graded ring, we have an augmentation ideal given by the projection onto the degree $0$ elements $\epsilon:A\to \ZZ$. Let $A^{\wedge}$ denote the completion at the augmentation ideal. 
%
%Let $k$ be a field and let $A$ be an $\NN$-graded $k$-algebra of finite type. Then $A$ is naturally augmented, with augmentation map $\epsilon:A\to k$ given by projecting onto the degree zero component $\epsilon(a)=a_0$. The completion at the augmentation ideal is again just $A$, so we can consider $A$ as a pro-algebra. The category of pro-modules over the pro-algebra $A$ is the limit ${\rm ProMod}_{A}:=\lim_{i\to \infty}{\rm Mod}_{A_i}$, where the structure maps ${\rm Mod}_{A_i}\to {\rm Mod}_{A_j}$ send $M_i\mapsto M_i\otimes _{A_i}A_j$. There is an obvious functor ${\rm Mod}_{A}\to {\rm ProMod}_{A}$ sending $M\mapsto \lim_{i\to \infty}M\otimes_{A}A_i$. Since $A_j$ is flat over $A$, this map is exact.

\begin{lemma}\label{Thefiberlem}
Let $A$ be a graded ring of finite type that is torsion free in each degree (as a $\ZZ$-module). We have a diagram in $\sm\Sp$
\begin{equation}\label{fiberdiag}
\xymatrix{
\mathcal{F}\ar[r] \ar[d]^-{\simeq} &\widehat{H}A(n)\ar[r]^-{\mathcal{R}}\ar[d]^-{\mathcal{I}} & H(Z\Omega^*(-,A[n]))\ar[d]
\\
\mathcal{G}\ar[r]&\Sigma^nH\underline{A}\ar[r] & H(\Omega^*(-,A[n])),
}
\end{equation}
where the square on the right is homotopy cartesian, both horizontal compositions are homotopy fiber sequences, the left vertical map is a weak equivalence, and the canonical map $\mathcal{G}\to \Sigma^nHA$ is weakly equivalent to the Bockstein map $\Sigma^{n-1}H(\underline{A\otimes_{\ZZ}\CC})\overset{\beta}{\to} \Sigma^nH\underline{A}$ . In particular, we have induced maps on cohomology
\begin{equation}\label{fiberdiag2}
\xymatrix{
&\widehat{H}A(n)\ar[r]^-{\mathcal{R}}\ar[d]^-{\mathcal{I}} & \Omega^n_{cl,bas}(-,A)\ar[d]^-{q}
\\
H^{n-1}(-,A\otimes_{\ZZ}\CC^{\times})\ar[r]^-{\beta}\ar[ru]^-{j} &H^n(-;A)\ar[r]^-{c} & H^n(-;A\otimes_{\ZZ} \CC),
}
\end{equation}
where the diagram is commutative and the two horizontal compositions are exact. The map $j$ is obtained by inverting the left vertical arrow in \eqref{fiberdiag} at the level of cohomology and composing with the top left horizontal arrow in \eqref{fiberdiag}. The map $c$ is the complexification map. The map $q$ is the map on cohomology induced by the inclusion 
$Z\Omega^*(-,A)\to \Omega^*(-,A)$.
\end{lemma}
\proof
Since $\sm\Sp$ presents a stable $\infty$-category, a square is homotopy cartesian if and only if the induced map on homotopy fibers of a pair of parallel arrows is a weak equivalence \cite[Proposition 1.1.3.4, Lemma 1.2.4.14]{LurieHA}. Let $\mathcal{F}$ and $\mathcal{G}$ denote the homotopy fibers of the top, respectively bottom, horizontal arrows in \eqref{fiberdiag}. Since the right square in \eqref{fiberdiag} is homotopy cartesian, the induced map $\mathcal{F}\to \mathcal{G}$ is a weak equivalence. The claim will follow provided we can show that $\beta$ is weakly equivalent to the canonical map out of the homotopy fiber. The homotopy fiber of the map $\Sigma^nH\underline{A}\to H(\Omega^*(-,A[n]))$ can be computed as follows. By the Poincar\'e lemma, the inclusion $\underline{A\otimes_{\ZZ}\CC}[n]\to \Omega^*(-,A[n])$ is a quasi-isomorphism of sheaves of complexes. Therefore the induced map $\Sigma^nH\underline{A}=H(\underline{A\otimes_{\ZZ}\CC}[n])\to H(\Omega^*(-,A[n]))$ is a weak equivalence in $\sm\Sp$. Hence the map $\Sigma^nH\underline{A}\to H(\Omega^*(-,A)[n])$ is weakly equivalent to the complexification map $c:\Sigma^nH\underline{A}\to H(\underline{A}\otimes_{\ZZ}\CC)$ and the homotopy fiber of the latter map produces the Bockstein \eqref{bockstien}.  
\endofproof

Next, we define the subring of integral elements in $I_{A}^*(\mathfrak{g}^{\vee})$. 
\begin{definition}
Let $A$ be a graded ring of finite type that is torsion free in each degree (as a $\ZZ$-module). We call an element in $I_{A}^*(\mathfrak{g}^{\vee})$ \emph{$H$-integral} if it is in the image of ring homomorphism
$$c:H^*(BG;A)\to H^*(BG;A\otimes_{\ZZ}\CC)\cong I_{A}^{*}(\mathfrak{g}^{\vee}),$$
where we have used Corollary \ref{chernweil} for identification on the right.
We denote the image of $c$ by $I_A^*(\mathfrak{g}^{\vee})_{H}\subset I_A^*(\mathfrak{g}^{\vee})$. 
\end{definition}

The integral elements in $I_A^*(\mathfrak{g}^{\vee})$ can be described geometrically via Chern--Weil theory. Given an element $p\in I_A^*(\mathfrak{g}^{\vee})$, a smooth manifold $M$, and a principal $G$-bundle with connection $M\to BG$, fix a connection $\mathcal{A}$ with curvature $\mathcal{F}$. By integrating over simplices, the differential form $p(\mathcal{F})$ corresponds to a singular cochain with values in $A\otimes_{\ZZ}\CC$. If this cochain has values in the graded ring $A\subset A\otimes _{\ZZ}\CC$ on cycles, then $p$ is integral.
%$$(\sigma:\Delta^n\to M)\mapsto \int_{\Delta^n}\sigma^*p(\mathcal{F}).$$
% where $\mathcal{F}$ denotes the curvature of the connection. Note that the power series $p(\mathcal{F})$ has only finitely many nonvanishing terms for degree reasons. Then $p$ is integral if the above integral is in the image of the canonical map $A\to A\otimes_{\ZZ}\CC$ whenever $\sigma$ is a cycle. The ring of integral elements can be computed using elementary representation theory.
% 
\begin{remark}\label{integralel}
The entire ring of integral elements can be described as follows. Let $A$ be a graded ring of finite type and assume that $A$ is torsion free in each degree (as a $\ZZ$-module). Let $G$ be compact and connected. Fix a maximal torus $T\subset G$ and let $\mathfrak{t}\subset \mathfrak{g}$ denote the corresponding Cartan algebra. Let $\lambda_1,\hdots,\lambda_n\in \mathfrak{t}^{\vee}$ be a basis of integral elements, i.e., $\lambda_i(H)\in \ZZ$ for all $H\in {\ker}({\rm exp})$. Let $W(G)$ denote the Weyl group. Then we have 
 $$I_A^*(\mathfrak{g})_{H}=\ZZ[[\lambda_1,\hdots,\lambda_n]]^{W(G)}\widehat{\otimes} \overline{A}, \quad |\lambda_i|=2	$$
Indeed, by Chern-Weil theory, we have an isomorphism 
$$H^*(BG;\CC)\cong I^{*}(\mathfrak{g}^{\vee})\cong \CC[\lambda_1,\hdots,\lambda_n]^{W(G)}.$$
 The complexification map $\ZZ[\lambda_1,\hdots,\lambda_n]^{W(G)}\cong H^*(BG;\ZZ)\to H^*(BG;\CC)$ is injective and can be identified with the canonical inclusion. The result then follows by completing the tensor product with $\overline{A}$. The resulting map is still injective, since $A$ is levelwise flat over $\ZZ$.
\end{remark}

Since the canonical inclusion $\iota:\underline{A}\to \Omega^*(-,A)$ is a map of presheaves of graded rings, the induced map of Eilenberg--MacLane spectra $H\underline{A}\to H(\Omega^*(-,A))$ is a map of ring spectra, which follows from \cite[Theorem 1.1]{Schipley}. The differential refinement 
$$\widehat{H}A=\bigvee_{n\in \ZZ}\widehat{H}A(n)$$
admits a compatible ring structure \cite[page 109]{Bunke}. In particular, for any smooth stack $X$, the curvature map $\mathcal{R}:\widehat{H}^*(X;A)\to \Omega^*_{cl}(X,A)$ is a map of graded rings. 
\begin{remark}
Note that some care must be taken with our chosen models for the Eilenberg--MacLane spectrum functor. More precisely, the Eilenberg--MacLane functor $\mathbb{H}$ considered in \cite{Schipley} is the composition of left and right adjoints in three Quillen equivalences. Our model of $H$ is only one of these functors. The precise statement is that the map $\mathbb{H}\underline{A}\to \mathbb{H}(\Omega^*(-,A))$ is a map of ring spectra. This technical point is not that important, since all these Quillen equivalences induce equivalences at the level of homotopy categories, and we are working at the level of cohomology anyways. More precisely, we are free to work with any model of $H$, keeping in mind that the ring structure on cohomology is transported along some isomorphism.
\end{remark}

We are now ready to prove Theorem \ref{theoremA}.

\medskip
{\it Proof of Theorem \ref{theoremA}\hspace{7pt}}
The curvature map (Definition \ref{diffcoh}) takes the form 
$$\mathcal{R}:\widehat{H}^n(\mathbf{B}_{\nabla}G;A)\to \Omega^n(\mathbf{B}_{\nabla}G,A)_{cl,bas}.$$
By Corollary \ref{basicfrmsbg}, we have an identification $\Omega^n_{cl,bas}(\mathbf{B}_{\nabla}G,A)\cong I^n_A(\mathfrak{g}^{\vee})$, so that $\mathcal{R}$ gives a ring homomorphism
$$\mathcal{R}:\widehat{H}^{*}(\mathbf{B}_{\nabla}G;A)\to I^*_A(\mathfrak{g}^{\vee}).$$
Next, we claim that $\mathcal{R}$ surjects onto integral elements. Clearly the image of $\mathcal{R}$ is contained in the subring $I^*_A(\mathfrak{g}^{\vee})_{H}$, by the commutativity $q\mathcal{R}=c\mathcal{I}$ in \eqref{fiberdiag2}. Fix an arbitrary element $p\in I^n_A(\mathfrak{g}^{\vee})_{H}$ of total degree $n$. Let $Q$ be a cofibrant replacement functor. Choose a derived map $x:Q\mathbf{B}_{\nabla}G\to K(\underline{A}[n])$ whose class in $H^n(BG;A)$ is an integral lift of $p$. Choose also a derived map $\tilde p:Q\mathbf{B}_{\nabla}G\to \Omega_{\rm cl}^n(-,A)$ whose class is the given class $p$. Then the composition 
$$Q\mathbf{B}_{\nabla}G\overset{x}{\to} K(\underline{A}[n])\overset{c}{\to} K(\underline{A}\otimes_{\ZZ}\CC[n])\overset{\iota}{\to} K(\Omega^*(-,A[n])) $$
 is homotopic to the composition 
 $$i\tilde p:Q\mathbf{B}_{\nabla}G\overset{\tilde p}{\to} \Omega_{\rm cl}^n(-,A)=K(Z\Omega^*(-,A[n]))\overset{i}{\to} K(\Omega^*(-,A[n])),$$
where $i$ is induced by the inclusion of cycles in degree zero and $\iota$ is induced by the inclusion of locally constant functions. That the maps are homotopic follows from the fact that $x$ is an integral lift of $p$.
A choice of such a homotopy induces, by the universal property of the homotopy pullbacks, a map $y:Q\mathbf{B}_{\nabla}G\to \mathbf{B}^n_{\nabla}A$ (See Lemma \ref{bnlemma} for notation) that makes the diagram commute up to homotopy. The class of $y$ satisfies $\mathcal{R}(y)=p$ by construction. This proves the assertion about $\mathcal{R}$ being surjective.

Next, consider the map $j:H^{*-1}(BG;A\otimes_{\ZZ}\CC^{\times})\to \widehat{H}^*(\mathbf{B}_{\nabla}G;A)$ in Lemma \ref{Thefiberlem}. By the exact sequence 
$$H^{n-1}(BG;A\otimes_{\ZZ}\CC^{\times})\overset{j}{\to} \widehat{H}^n(\mathbf{B}_{\nabla}G;A)\overset{\mathcal{R}}{\to} I_A^*(\mathfrak{g}^{\vee}),$$
the kernel of $\mathcal{R}$ is the image of $j$. By the commutativity $\mathcal{I}j=\beta$ in Lemma \ref{Thefiberlem}, it follows that $\mathcal{I}$ maps the image of $j$ isomorphically onto the image of $\beta$. The image of the Bockstein homomorphism
$$\beta:H^{n-1}(BG;A\otimes_{\ZZ}\CC^{\times})\to H^{n}(BG;A)$$
is precisely the $\ZZ$-torsion in $H^{n}(BG;A)$. 

%Finally, when $G$ is compact and connected the sequence splits since $I_A^*(\mathfrak{g})_H$ is levelwise projective as a pro $A$-module. Indeed, by Remark \ref{integralel}, we have an inclusion
%%\begin{equation}\label{inclinv} I^*_A(\mathfrak{g})_H=\ZZ[[\lambda_1\hdots,\lambda_n]]^{W(G)}\widehat{\otimes}A \into \ZZ[[\lambda_1,\hdots,\lambda_n]]\widehat{\otimes}A.
%%\end{equation} 
%%The canonical inclusion of $W(G)$-invariants
%%$$\ZZ[[\lambda_1,\hdots,\lambda_n]]^{W(G)}\into \ZZ[[\lambda_1,\hdots,\lambda_n]]$$
%can be identified with the procompletion of the correpsonding map on polynomial rings (\cite{Atiyah}). The inclusion is split in the category of prorings by the limit of the map $a_n(x)=\frac{1}{\vert W(G)\vert}\sum_{\sigma\in W(G)}\sigma x, x\in \CC[\lambda_1,\hdots,\lambda_n]/I(G)^n$ \footnote{Note that the procompletions of $R(T)=\CC[\lambda_1,\hdots,\lambda_n]$ at $I(G)$ and $I(T)$ coincide (\cite{atiyha})}. The map $a=\lim_{n\to \infty}a_n$ The induced map $a\widehat{\otimes}A$ splits the inclusion \eqref{inclinv} in the category of pro $A$-modules.

\endofproof

To illustrate what can happen in the noncompact case, we have the following example. 
\begin{example}
Let $A=\ZZ$ and let $G=\ZZ$. Then $\mathbf{B}_{\nabla}\ZZ=\mathbf{B}\ZZ$. We have $H^*(\mathbf{B}\ZZ;\ZZ)\cong H^*(S^1;\ZZ)\cong \Lambda_{\ZZ}(x)$, $|x|=1$, $H^*(\mathbf{B}\ZZ;\RR)\cong \Lambda_{\RR}(x)$, $|x|=1$, $I^*(\mathfrak{g}^{\vee})\cong 0$ and $H^*(\mathbf{B}\ZZ;\CC^{\times})\cong \CC^{\times}[0]\oplus \CC^{\times}[1]$. It follows from Theorem \ref{theoremA} and Lemma \ref{Thefiberlem} that the map $\mathcal{I}:\widehat{H}(\mathbf{B}\ZZ;\ZZ)\to H^*(\mathbf{B}\ZZ;\ZZ)$ is not surjective and ${\rm im}(\mathcal{I})={\rm im}(\beta)$. From the Bockstein sequence, we learn that the map $\beta=0$ in all cases, so that $\mathcal{I}=0$. Therefore $j:H^{*-1}(\mathbf{B}\ZZ;\CC^{\times})\to \widehat{H}^*(\mathbf{B}\ZZ;\ZZ)$ is an isomorphism, so that $\widehat{H}^*(\mathbf{B}_{\nabla}\ZZ;\ZZ)\cong \CC^{\times}[1]\oplus \CC^{\times}[2]$. In this case, the ring structure is trivial.
\end{example}

\section{The differential $K$-theory of $\mathbf{B}_{\nabla}G$}
\def\Perm{\mathscr{P}{\rm erm}}
\def\Gammas{\Gamma\text{-}\sset}

In this section, we prove the second main result. First, we recall how to define the various differential $K$-theory spectra on smooth stacks. We recall from \cite{Elmendorf} the functor $\mathcal{K}:{\rm Perm}\to \Sp$ from the multicategory of permutative categories to the multicategory of symmetric spectra. The underlying functor of this multifunctor is equivalent to Segal's construction of a spectrum from a permutative category \cite[Theorem 4.6]{Elmendorf}. The following constructions are due to \cite{BNV} in the setting of quasi-categories. Here we briefly explain the analogous constructions using model categories.

Let $\mathscr{P}{\rm erm}$ denote the category of permutative categories. We define a \emph{smooth permuative category} to be a functor  $\Cart^{\rm op}\to \mathscr{P}{\rm erm}$.

\begin{definition}
Let $P\in \sm\Perm$ be a smooth permutative category. We define its $K$-theory as the smooth spectrum by applying the functor $\mathcal{K}$ objectwise:
$$
U\mapsto \mathcal{K}P(U), \quad U\in \Cart.
$$
This produces a functor, which we denote by the same symbol:
$$\mathcal{K}:\sm\Perm\to \sm\Sp.$$
\end{definition}

\begin{example}
The following two smooth permutative groupoids will be used in the construction of $K$-theory. Consider the smooth permutative groupoid of vector bundles ${\rm Vect}^{\oplus}$, which sends $U\in \Cart$ to the permutative groupoid 
$${\rm Vect}^{\oplus}(U)=\coprod_{n\geq 0}\mathbf{B}C^{\infty}(U,{\rm GL}_n(\CC)).$$
The monoidal operations given by the usual direct sum of matrices. Since pullback smooth maps to matrices commutes with direct sums, this indeed defines a smooth permutative category (see \cite[Example, page 15]{Elmendorf}).

We also have the smooth permuative groupoid of vector bundles with connection ${\rm Vect}^{\oplus}_{\nabla}$, which sends $U\in \Cart$ to the groupoid 
$${\rm Vect}^{\oplus}_{\nabla}(U)=\coprod_{n\geq 0}\Omega^1(U,\mathfrak{gl}_n)\rtimes C^{\infty}(U,{\rm GL}_n),$$
where $\Omega^1(U,\mathfrak{gl}_n)\rtimes C^{\infty}(U,{\rm GL}_n)$ denotes the action groupoid of smooth functions acting on 1-forms by gauge transformations. Both ${\rm Vect}^{\oplus}$ and ${\rm Vect}^{\oplus}_{\nabla}$ admits the structure of a bipermutative (or rig) structure under the tensor product $\otimes$. In the later case, the multiplicative structure is given on objects by 
$$\mathcal{A}\otimes \mathcal{B}(e_i\otimes e_j)=\mathcal{A}e_i\otimes e_j+e_i\otimes \mathcal{B}e_j, \quad(g\otimes h)(e_i\otimes e_j)=g(e_i)\otimes h(e_j)$$
for $\mathcal{A}\in \Omega^1(U,\mathfrak{gl}_n),\mathcal{B}\in \Omega^1(U,\mathfrak{gl}_m),g\in C^{\infty}(U,GL_n), h\in C^{\infty}(U,GL_m)$, and $e_i\in \RR^n,e_j\in \RR^m$ are the standard basis vectors. The distributive maps are the obvious ones, induced by the isomorphism $\RR^n\otimes (\RR^m\oplus \RR^{\ell})\to (\RR^n\otimes \RR^m)\oplus (\RR^n\otimes \RR^{\ell})$. By \cite[Theorem 3.7, Corollary 3.9]{Elmendorf}, it follows that the smooth spectra obtained by applying the $K$-theory functor to the above bipermutative categories admit an $E_{\infty}$-ring structure. 
\end{example}

%Composing the smooth permutative groupoid ${\sf Vect}^{\oplus}$ with Segal's nerve functor $N_{\Gamma}:{\rm Perm}\to \Gamma\mathscr{S}$ produces a smooth $\Gamma$-space. This allows us to make the following definition
%The upshot of this abstract nonsense is that it allows us to define a groupoid of vector bundles on any smooth stack $X\in \smsset$. Indeed, since ${\rm Vect}^{\times}$ is 1-truncated, being the nerve of a groupoid, the derived mapping simplicial set 
%$\mathbb{R}\map(X,{\rm Vect}^{\times})$
%is also 1-truncated. It follows at once that $\mathbb{R}\map(X,{\rm Vect}^{\times})$ is equivalent to the nerve of its fundamental groupoid. Since ${\rm Vect}^{\times}$ is fibrant in the local projective model structure, we have the following definition

\begin{definition}
Let $X\in \sm\mathscr{S}$ and let $Q$ be a cofibrant replacement functor.
\begin{enumerate}
%\item We define the $\Gamma${\bf -space of vector bundles on} $X$ as 
%$${\rm Vect}^{\oplus}(X):=\Gamma\mathscr{S}(X,N_{\Gamma}({\rm Vect}^{\oplus}))$$
%A map $\langle 1\rangle\to {\rm Vect}^{\oplus}(X)$ is called a vector bundle on $X$. 
%\item We define the $\Gamma${\bf -space of vector bundles with connection on} $X$ as 
%$${\rm Vect}^{\oplus}(X):=\Gamma\mathscr{S}(X,N_{\Gamma}({\rm Vect}_{\nabla}^{\oplus}))$$
\item We define the smooth connective $K$-theory of $X$ as 
$$k(X):=\mathcal{K}{\rm Vect}^{\oplus}(QX),$$
\item We define the differential connective $K$-theory of $X$ as 
$$k_{\nabla}(X):=\mathcal{K}{\rm Vect}_{\nabla}^{\oplus}(QX).$$
\end{enumerate}
\end{definition}

\def\smsp{{\rm C}^{\infty}\mathscr{S}{\rm p}}

Next, we turn our attention to the Hopkins--Singer differential $K$-theory. We recall the definition by \cite{BNV}. Let $\Omega^*(-,\CC[u])$, $|u|=2$, be the complex of even and odd forms. Taking cycles in degree zero and applying the Eilenberg-MacLane spectrum functor $H$ yields a smooth spectrum $HZ\Omega^*(-,\CC[u])$, which serves as the differential form component of differential (connective) $K$-theory. The inclusion of constant forms $\underline{\CC}[u]\into \Omega^*(-,\CC[u])$ is a quasi-isomorphism of sheaves of differentially graded commutative algebras, by the Poincar\'e lemma, hence induces an equivalence of smooth spectra $H\underline{\CC}[u]\to H\Omega^*(-,\CC[u])$. 
%Taking the relative category of smooth symmetric spectra with local stable weak equivalences and applying hammock localization, we obtain a map the other way
%\begin{equation}\label{bunkediffer}
%\iota:H\tau_{\leq 0}\Omega^*[u]\to H\CC[u]
%\end{equation}
%which corresponds to the hammock 
%$$H\tau_{\geq 0}\Omega^*[u]\leftarrow H\CC[u]\overset{=}\to H\CC[u].$$

The Chern character form (defined via Chern-Weil theory) defines a natural transformation ${\rm ch}:{\rm Vect}_{\nabla}\to Z\Omega^*(-,\CC[u])$. Regarding the smooth set on the right as the discrete bipermutative category with $\oplus$ given by the sum of differential forms and $\otimes$ given by the wedge product, this defines a strong monoidal functor, hence induces a morphism of smooth $E_{\infty}$-ring spectra
$$k_{\nabla}\overset{\rm ch}{\to} HZ\Omega^*(-,\CC[u]).$$
%with the latter morphism is induced by quasi-isomorphism $Z(\Omega^*(-,\CC[u]))\to \tau_{\leq 0}\Omega^*(-,\CC[u])$. The $E_{\infty}$-structure on $H(\tau_{\geq 0}\Omega^*(-,\CC[u]))$ is the obvious one, induced by the wedge product of forms. 
\begin{lemma}\label{chernchar}
The Chern character map ${\rm ch}:k_{\nabla}\to HZ\Omega^*(-,\CC[u])$ induces a map of $E_{\infty}$-ring spectra
$$
c:=|{\rm ch}|:\vert k_{\nabla}\vert \to |HZ\Omega^*(-,\CC[u])|,
$$
where $|k_{\nabla}|\simeq k$ and $
|HZ\Omega^*(-,\CC[u])|\simeq H\CC[u]$.
\end{lemma}
\proof
That the map induces a map  of $E_{\infty}$-ring spectra upon application of $|\cdot |$ follows from Proposition \ref{commeinfring} (the ring structure is provided by Corollary \ref{unitshapes}). It remains to verify the second statement about the identification of the shapes. To see this, apply $\vert \cdot \vert$ to the weak equivalence $HZ\Omega^*(-,\CC[u])\to H\tau_{\geq 0}\Omega^*(-,\CC[u])$. The result then follows by \cite[Lemma 4.4, part 5]{BNV}, taking $m=0$ and $C=\CC[u]$ there, and \cite[Corollary 6.5]{BNV}. 

 \endofproof
  
The definition of (connective) differential $K$-theory according to \cite{BNV} is the homotopy pullback 
$$
\xymatrix{
\widehat{k}\ar[r]^-{\mathcal{I}}\ar[d]_-{\mathcal{R}} & H\tau_{\leq 0}\Omega^*(-,\CC[u])\ar[d]^-{\iota}
\\
k\ar[r]^-{c} & H\Omega^*(-,\CC[u]).
}
$$
In order to simplify some of the calculations, we will work with the following weakly equivalent model.

\begin{definition}
The Hopkins--Singer (connective) differential $K$-theory in degree $n$ is obtained from the data $(k,\CC[u], c)$, via the homotopy pullback in the category of sheaves of spectra. 
$$
\xymatrix{
\widehat{k}(n)\ar[r]^-{\mathcal{I}}\ar[d]_-{\mathcal{R}} & H(Z\Omega^*(-,\CC[u][n]))\ar[d]^-{\iota}
\\
\Sigma^nk\ar[r]^-{\Sigma^nc} & H\Omega^*(-,\CC[u][n])
}
$$
where $c$ is the Chern character map in Lemma \ref{chernchar} and $\iota$ is induced by the canonical map $Z\Omega^*(-,\CC[u])\to \Omega^*(-,\CC[u])$.
\end{definition}

\begin{remark}
Since both $c$ and $\iota$ are maps of $E_{\infty}$-ring spectra, $\widehat{k}$ is a multiplicative extension of $k$ (in the sense of \cite{Bunke}). At the level of cohomology, the multiplicative structure is compatible with topological $K$-theory and differential forms: that is, the maps $\mathcal{I}$ and $\mathcal{R}$ induce homomorphisms of graded rings at the level of cohomology. This is spelled out in \cite{Bunke}. See also \cite{GradySati} in the case of $KO$-theory.
\end{remark}

Recall that the exponential sequence 
$$0\to \ZZ\to \CC\to \CC^{\times}\to 0$$
induces a long fiber/cofiber sequence, by smashing with the corresponding Moore spectra:
$$\hdots \to\Sigma^{-1}k_{\CC^{\times}} \overset{\beta}{\to} k\to k_{\CC}\to k_{\CC^{\times}}\overset{\beta}{\to} \Sigma k \hdots$$
As in the case of cohomology, we have the following. The proof is completely analogous to Lemma \ref{Thefiberlem} and we will omit it.
\begin{lemma}\label{Thefiberlem2}
We have a diagram in $\sm\Sp$
\begin{equation}\label{fiberdiag3}
\xymatrix{
\mathcal{F}\ar[r] \ar[d]^-{\simeq} &\widehat{k}(n)\ar[r]^-{\mathcal{R}}\ar[d]^-{\mathcal{I}} & H(Z\Omega^*(-,\CC[u][n]))\ar[d]
\\
\mathcal{G}\ar[r]&\Sigma^n k\ar[r]^-{\Sigma^nc} & H(\Omega^*(-,\CC[u][n]))
}
\end{equation}
where the square on the right is homotopy cartesian, both horizontal compositions are fiber sequences, the left vertical map is a weak equivalence, and the canonical map $\mathcal{G}\to k$ is weakly equivalent to the Bockstein map $\Sigma^{n-1}k_{\CC^{\times}}\overset{\beta}{\to} \Sigma^nk$. In particular, we have induced maps on cohomology
\begin{equation}\label{fiberdiag4}
\xymatrix{
&\widehat{k}^n\ar[r]^-{\mathcal{R}}\ar[d]^-{\mathcal{I}} & \Omega^n_{cl,bas}(-,\CC[u])\ar[d]^-{q}
\\
k^{n-1}_{\CC^{\times}}\ar[r]^-{\beta}\ar[ru]^-{j} &k^n\ar[r]^-{c} & H^n(-;\CC[u])
}
\end{equation}
where the diagram is commutative and the two horizontal compositions are exact. The map $j$ is obtained by inverting the left vertical arrow in \eqref{fiberdiag3} at the level of cohomology and composing with the top left horizontal arrow in \eqref{fiberdiag3}. The map $c$ is the complexification map. The map $q$ is the map on cohomology induced by the inclusion 
$Z\Omega^*(-,\CC[u])\to \Omega^*(-,\CC[u])$.
\end{lemma}
%Next, we will define the ring of $k$-integral invariant polynomials. Let $A=\ZZ[u]$, $|u|=2$ be the graded ring of polynomials in $u$ with integer coefficients. Then the pro-completed tensor product of graded rings $I^*(\mathfrak{g}^{\vee})\widehat{\otimes}_{\ZZ}\ZZ[u]$ can be identified with the graded ring $I(\mathfrak{g}^{\vee})^{\wedge}[v]$, $|v|=2$, where $I(\mathfrak{g}^{\vee})^{\wedge}$ is the pro-completion 
%$$I(\mathfrak{g}^{\vee})^{\wedge}=\varprojlim_k I^{\leq k}(\mathfrak{g}^{\vee}),$$ 
%taken in the category of rings (forgetting the grading). The variable $v$ keeps track of the total grading in the tensor product.

\begin{definition}
By Corollary \ref{chernweil}, we have an isomorphism 
$$H^*(BG;\CC[u])\cong I^*_{\CC[u]}(\mathfrak{g}^{\vee}).$$
The ring on the right is described explicitly in Example \ref{invariantpolys}. The Chern character map $c$ gives a homomorphism of graded rings
$$c:k^*(BG)\to I^*_{\CC[u]}(\mathfrak{g}^{\vee}).$$
We denote by $I^*_{\CC[u]}(\mathfrak{g}^{\vee})_{k}$ the subring obtained by taking the image of the above homomorphism. We call an element $f\in I^*_{\CC[u]}(\mathfrak{g}^{\vee})_{k}$ a \emph{$k$-integral element}.  
\end{definition}

We are now ready to prove Theorem \ref{theoremB}.
\medskip

{\it Proof of Theorem \ref{theoremB}\hspace{7pt}}
By the long homotopy fiber/cofiber sequence associated to the sequence of sheaves of spectra $\Sigma^{n-1}k_{\CC^{\times}}\to \hat k(n)\to H(Z\Omega^*(-,\CC[u][n])$, we see that $\hat k^n(\mathbf{B}_{\nabla}G)$ sits in an exact sequence
\begin{equation}\label{lesdiffk}
k^{n-1}_{{\CC}}(\mathbf{B}_{\nabla}G)\to k^{n-1}_{\CC^{\times}}(\mathbf{B}_{\nabla}G)\to \hat k^n(\mathbf{B}_{\nabla}G)\to H^0(Z\Omega^*(\mathbf{B}_{\nabla}G,\CC[u][n]))
\end{equation}
Since $k_{\CC^{\times}}$ is a homotopy invariant sheaf of spectra, it follows from Proposition \ref{stablecoh} and Example \ref{clssbg} that $k^*_{\CC^{\times}}(\mathbf{B}_{\nabla}G)\cong k^*_{\CC^{\times}}(BG)$. Similarly, $k^*_{\CC}(\mathbf{B}_{\nabla}G)\cong k^*_{\CC}(BG)$.
 Hence, the exact sequence \eqref{lesdiffk} reduces to an exact sequence 
$$
k^{n-1}_{\CC^{\times}}(BG) \to \widehat{k}^n(\mathbf{B}_{\nabla}G)\overset{\mathcal{R}}{\to} H^0(Z\Omega^*(Q(\mathbf{B}_{\nabla}G),\CC[u][n])),
$$
 By Corollary \ref{basicfrmsbg}, we have an isomorphism
$$H^0(Z\Omega^*(\mathbf{B}_{\nabla}G,\CC[u]))=\Omega^n(\mathbf{B}_{\nabla}G,\CC[u])_{cl,bas}\cong I_{\CC[u]}^n(\mathfrak{g}^{\vee}).$$
Hence, the sequence \eqref{lesdiffk} reduces to an exact sequence
\begin{equation}\label{exactcwev}
 k^{n-1}_{\CC^{\times}}(BG) \overset{j}{\to} \widehat{k}^n(\mathbf{B}_{\nabla}G)\overset{\mathcal{R}}{\to} I^n_{\CC[u]}(\mathfrak{g}^{\vee}).
\end{equation}
By the commutativity $\mathcal{I}j=\beta$ in Lemma \ref{fiberdiag2}, it follows that $\mathcal{I}:\widehat{k}^n(\mathbf{B}_{\nabla}G)\to k^n(BG)$ maps the kernel of $\mathcal{R}$ onto the image of $\beta$.

Next, we show that $\mathcal{R}$ surjects onto $I_{\CC[u]}^*(\mathfrak{g}^{\vee})_{k}$. By definition, any element $p\in I_{\CC[u]}^n(\mathfrak{g}^{\vee})_{k}$ can be lifted through the Chern character to a cocycle in $k$-theory 
\begin{equation}\label{intliftk}
x:\Sigma_+^{\infty}Q(\mathbf{B}_{\nabla}G)\to \Sigma_+^{\infty}Q\delta(BG) \to \Sigma^nk.
\end{equation}
Here we have used the homotopy adjunction $\delta\dashv |\cdot |$ in Proposition \ref{stablecoh}, Example \ref{clssbg}, and homotopy invariance of $\Sigma^nk$ to conclude that the map $\Sigma_+^{\infty}Q(\mathbf{B}_{\nabla}G)\to \Sigma_+^{\infty}Q\delta(BG)$ induces a weak equivalence when mapping into $\Sigma^nk$, so the lift $x$ can be chosen to be of the form \eqref{intliftk}.
Moreover, $p$ can be lifted to a cocycle in de Rham cohomology $\tilde p:\Sigma^{\infty}_+Q(\mathbf{B}_{\nabla}G)\to HZ(\Omega^*(-,\CC[u]))$, which follows from Corollary \ref{basicfrmsbg}. The universal property of the homotopy pullback produces a cocycle in $\hat k(n)$ that maps to the given elements. Hence, $\mathcal{R}$ is surjective onto $I^n_{\CC[u]}(\mathfrak{g}^{\vee})_k$, as claimed.

It remains to show that $\mathcal{R}$ is an isomorphism in degrees $n\leq 0$ when $G$ is compact and connected. By Chern--Weil theory, the group $k^{*}_{\CC}(BG)\cong H^*(BG,\CC[u])$ vanishes in odd degrees. Hence, for $n\leq 0$ even, we have a short exact sequence
\begin{equation}\label{exactcwev}
0\to k^{n-1}_{\CC^{\times}}(BG) \overset{j}{\to} \widehat{k}^n(\mathbf{B}_{\nabla}G)\overset{\mathcal{R}}{\to} I^n_{\CC[u]}(\mathfrak{g}^{\vee})_k\to 0,
\end{equation}
%while for $n$ odd, we have the exact sequence
%\begin{equation}\label{exactcwodd}
%k^{n-1}(BG)\to k^{n-1}_{\CC}(BG)\to k^{n-1}_{\CC^{\times}}(BG) \overset{j}{\to} \widehat{k}^n(\mathbf{B}_{\nabla}G)\to 0.
%\end{equation}
By the universal coefficient theorem, we have an exact sequence 
\begin{equation}\label{exactuct}
0\to k^{n-1}(BG)\otimes \CC^{\times}\to k^{n-1}_{\CC^{\times}}(BG)\to {\rm Tor}(k^{n}(BG),\CC^{\times})\to 0.
\end{equation}
By \cite[Proposition 2.6]{Greenlees}, $k^n(BG)=0$ when $n\leq 0$, $n$ odd, and $k^n(BG)=R(G)^{\wedge}$ when $n\leq 0$, $n$ even. Since $R(G)^{\wedge}$ is torsion free, it follows that for $n$ even, $n\leq 0$, we have $k^{n-1}(BG)\otimes \CC^{\times}\cong k^{n-1}_{\CC^{\times}}(BG)\cong 0$. Hence, $\mathcal{R}$ is an isomorphism in degrees $n\leq 0$, $n$ even.

%
%The entire argument works in the case where $G$ is not compact or connected, except for injectivity of the map $k^{n-1}_{\CC^{\times}}(BG)\to \widehat{k}^n(\mathbf{B}_{\nabla}G)$, since the odd degree cohomology groups of $H^*(BG,\CC[u])$ need not vanish. This proves the claim in the noncompact case.
\endofproof

We conclude this section by computing $\widehat{k}^0(\mathbf{B}_{\nabla}G)$ in the case where $G$ is compact and connected, proving Corollary \ref{corC}.

\medskip

{\it Proof of Corollary \ref{corC}\hspace{7pt}}
The commutativity of the diagram follows from the identification of $\mathcal{R}$ in Theorem \ref{theoremB} and the Atiyah--Segal completion theorem. Let $T\subset G$ be a maximal torus with corresponding Cartan algebra $\mathfrak{t}$. Let $\{\lambda_i\}$ be a basis of $\mathfrak{t}^{\wedge}$ consisting of integral elements. According to \cite[Theorem, page 19]{Atiyah}, we have an identification 
$$R(G)^{\wedge}=\ZZ[[\lambda_1,\hdots,\lambda_n]]^{W(G)}.$$
Under this identification, the Chern character map can be identified with the injective map 
$$
\ZZ[[\lambda_1,\hdots,\lambda_n]]^{W(G)}\to \CC[[\lambda_1,\hdots,\lambda_n]]^{W(G)}\cong H^0(BG;\CC[u]), \quad \lambda	_i\mapsto e^{\lambda_i}-1.
$$
On the other hand, as we have seen in Example \ref{invariantpolys}, we have isomorphisms
$$\CC[[\lambda_1,\hdots,\lambda_n]]^{W(G)}\cong I^0_{\CC[u]}(\mathfrak{g}^{\wedge})\cong  H^0(BG;\CC[u]),$$ 
where the second isomorphism is obtained from Corollary \ref{chernweil}. Under the above identifications, the Chern character map gives an isomorphism ${\rm ch}:R(G)^{\wedge}\to I^0_{\CC[u]}(\mathfrak{g}^{\wedge})_k$ onto the $k$-integral elements. Since $G$ is compact and connected, the result then follows immediately from Theorem \ref{theoremB}, taking $n=0$.

\endofproof

\appendix
\section{Sheaves of spectra}

The purpose of this appendix is to review known constructions which are used throughout this work and to prove some general results about sheaves of spectra, which are needed in the present work.

We will begin with some background. There are several adjoint Quillen functors relating $\sm\mathscr{S}$ and $\mathscr{S}$. The most primative of these functors is the constant sheaf functor. This functor admits both a left and right adjoint by Kan extension (in this case, they are just the colimit and limit functors, respectively). In fact, there is even a further right adjoint of the limit functor. Altogether, these four functors give $\sm\mathscr{S}$ the structure of a \emph{cohesive $\infty$-topos}, which was introduced in \cite{Urs}. 
\begin{proposition}\label{cohesive}
The $\infty$-topos  $\sm\mathscr{S}$ is \emph{cohesive} over $\mathscr{S}$. In fact, we have a quadruple Quillen adjunction at the level of the local projective model structure $\vert \cdot \vert\dashv \delta\dashv\Gamma\dashv \delta^{\dagger}$:
$$
\xymatrix{
\sm\mathscr{S}\ar@{{x}->}@<.4cm>[rr]^-{\vert\cdot \vert}\ar@<-.15cm>[rr]|-{\Gamma} && \mathscr{S}\ar@{_{(}->}@<-.15cm>[ll]|-{\delta} \ar@{_{(}->}@<.4cm>[ll]^-{\delta^{\dagger}}
}
$$
where the functors are defined by the following formulas:
\begin{itemize}
\item $\Gamma(X)=X(\RR^0)$
\item $\vert X\vert=\colim_{U\in \Cart^{\rm op}}X(U)$
\item $\delta(Y)(U)=Y$
\item $\delta^{\dagger}(Y)(U)=Y^{C^{\infty}(\RR^0,U)}$,
\end{itemize}
the functor $|\cdot |$ preserves finite products, $\delta$ and $\delta^{\dagger}$ are (homotopically) fully faithful. 
\end{proposition}
\proof
This is \cite[Proposition 3.4.18]{Urs}.
\endofproof

\begin{remark}
We will only use the derived functor of the colimit functor $|\cdot |$ (i.e., the homotopy colimit) in the present work. To avoid carrying around the notation for left derived functor, we will simply use the notation $|\cdot|$ for the derived functor throughout. 
\end{remark}

The following example is fundamental. 
\begin{example}\label{clssbg}
Let $G$ be a Lie group. Let $\mathbf{B}_{\nabla}G\in \sm\mathscr{S}$ denote the moduli stack of principal $G$-bundles, presented by the simplicial presheaf  $\mathbf{B}_{\nabla}G(U)=N(\Omega^1(U;\mathfrak{g}))\rtimes C^{\infty}(U,G))$: the nerve of the action groupoid where smooth functions act on 1-forms by gauge transformations. Then we have a weak equivalence
$$\vert \mathbf{B}_{\nabla}G\vert\overset{\simeq}{\to} BG,$$
where $BG$ is classifying space of $G$ (regarded as a simplicial set via the nerve). Similarly, the moduli stack of smooth principal $G$-bundles $\mathbf{B}G(U)=N(\ast\rtimes C^{\infty}(U,G))$ has shape $|\mathbf{B}G\vert \simeq BG$. 
\end{example}

Proposition \ref{cohesive} extends to smooth spectra in an obvious way. That is we again have the colimit functor $\vert \cdot \vert:\sm\Sp\to \Sp$ and the constant functor $\delta:\Sp\to \sm\Sp$, which fit into a Quillen adjunction $\vert \cdot \vert\dashv \delta$. In fact, the entire quadruple adjunction extends an $\infty$-adjunction at the level of quasi-categories. 
\begin{proposition}\label{stablecoh}
The stable $\infty$-category $\sm\Sp$ is \emph{cohesive} over $\Sp$. In fact, we have a quadruple $\infty$-adjunction $\vert \cdot \vert\dashv \delta\dashv\Gamma\dashv \delta^{\dagger}$:
$$
\xymatrix{
\sm\Sp\ar@{{x}->}@<.4cm>[rr]^-{\vert\cdot \vert}\ar@<-.15cm>[rr]|-{\Gamma} && \Sp\ar@{_{(}->}@<-.15cm>[ll]|-{\delta} \ar@{_{(}->}@<.4cm>[ll]^-{\delta^{\dagger}}
}
$$
where the functors are defined by the following formulas:
\begin{itemize}
\item $\Gamma(X)=X(\RR^0)$
\item $\vert X\vert=\colim_{U\in \Cart^{\rm op}}X(U)$
\item $\delta(Y)(U)=Y$
\item $\delta^{\dagger}(Y)(U)=Y^{C^{\infty}(\RR^0,U)}$,
\end{itemize}
the functor $|\cdot |$ preserves finite products, $\delta$ and $\delta^{\dagger}$ are fully faithful. 
\end{proposition}
\proof
 This is \cite[Proposition 2.6]{BNV}.
\endofproof
 The first goal in this appendix is to show that the adjunction $|\cdot |\dashv \delta$ induces an adjunction on $E_{\infty}$-algebras. This is essentially a direct application of \cite[Proposition 7.9]{Pavlov}. Let $\Sp_{E_{\infty}}$ denote the $\infty$-category of $E_{\infty}$-ring spectra, modeled as follows. Let $\Sigma_*$ be the associative operad in sets, whose $k$-th component is the symmetric group $\Sigma_k$ and whose multiproduct is defined by 
 $$m:\Sigma_k\times \Sigma_{j_1}\times \hdots \times \Sigma_{j_k}\to \Sigma_{j_1+\hdots +j_k}, \quad m(\sigma,\tau_{j_1},\hdots,\tau_{j_k})(x)=\sigma\tau_{j_i}(x),$$
 where $x\in \{j_1+\hdots +j_{i-1}, \hdots, j_1+\hdots +j_{i}\}$ and $\sigma\in \Sigma_k$ acts on the set $\{1,\hdots, j_1+\hdots+j_k\}$ by permuting the $k$ blocks specified by the $j_i$s. Let $E:\set\to \mathscr{C}{\rm at}$ be the functor that sends a set $X$ to the category with objects $X$ and exactly 1 morphism between any pair of objects. Let $E\Sigma_*$  be the corresponding $E_{\infty}$-operad in the category of small categories. The category of $E_{\infty}$-ring spectra can be modeled as the category of algebras over $E\Sigma_*$ (as described in \cite{Elmendorf}), which we denote by $\Sp_{E_{\infty}}$. This category admits a simplicial model structure, whose weak equivalences are objectwise stable equivalences and whose fibrations are objectwise positive stable fibrations \cite{Pavlov}.  Let $\sm\Sp_{E_{\infty}}$ denote the category of presheaves on $\Cart$ with values in $\Sp_{E_{\infty}}$, equipped with the projective model structure. 
 
We have the following

\begin{proposition}\label{commeinfring}
Let $\vert \cdot \vert:\sm\Sp_{E_{\infty}}\to \Sp_{E_{\infty}}$ denote the homotopy colimit functor. Let $u:\Sp_{E_{\infty}}\to \Sp$ and $u:\sm\Sp_{E_{\infty}}\to \sm\Sp$ denote the corresponding forgetful functors. Then we have a zig-zag of weak equivalence of spectra
$$u\vert \mathcal{E}\vert\simeq \vert u\mathcal{E} \vert$$
\end{proposition}
\proof
Apply \cite[Proposition 7.9]{Pavlov} to the case where $\mathscr{C}=\Sp$, $\mathcal{O}=E\Sigma_*$ and use the fact that $\Cart$ is sifted.
\endofproof

We say that a smooth spectrum $E\in \sm\Sp$ admits the structure of a smooth $E_{\infty}$-ring spectrum if it is in the essential image of the forgetful functor $\sm\Sp_{E_{\infty}}\to \sm\Sp$. Similarly, we say that a spectrum $E\in \Sp$ admits the structure of an $E_{\infty}$-ring spectrum if it is in the essential image of $\Sp_{E_{\infty}}\to \Sp$. We say that a map of (smooth) spectra $f:E\to F$ is a morphism of $E_{\infty}$-ring spectra if it lifts through the forgetful functor to a morphism of corresponding (smooth) $E_{\infty}$-ring spectra. As an immediate corollary we have the following.

\begin{corollary}\label{unitshapes}
Let $E\in \sm\Sp$ be a  smooth spectrum admitting the structure of a smooth $E_{\infty}$-ring spectrum. Then $\vert E\vert$ admits the structure of a $E_{\infty}$-ring spectrum such that the unit of the adjunction $\eta_{E}:E\to \delta\vert E\vert$ is a morphism of $E_{\infty}$-ring spectra. In particular, the unit induces an homorphism of rings
$$\eta_{E}:E(X)\to \delta\vert E\vert^*(X).$$ 
Moreover, if $E$ is homotopy invariant, then $\eta_{E}$ is a weak equivalence of $E_{\infty}$-ring spectra and $\eta_{E}$ is an isomorphism.
\end{corollary}
\proof
For the first claim, let $\mathcal{E}$ be the lift of $E$ to an $E_{\infty}$-ring spectrum. Consider the $E_{\infty}$-ring spectrum $\vert \mathcal{E}\vert$. By the commutativity of the diagram in Proposition \ref{commeinfring}, $|\mathcal{E}\vert$ defines an $E_{\infty}$-ring structure on $\vert E\vert$. The Quillen adjunction $|\cdot|\dashv \delta$ continuous to hold at the level $E_{\infty}$-ring spectra, since the functor $\delta$ manifestly preserves objectwise fibrations and acyclic fibrations, hence is right Quillen. Applying the forgetful functor to the $\mathcal{E}$-component of the unit of the adjunction recovers the unit $\eta_{E}$. Hence, $\eta_{E}$ is a morphism of $E_{\infty}$-ring spectra. This proves the first claim. If $E$ is homotopy invariant, then let $\mathcal{E}$ be a lift to an $E_{\infty}$-ring structure. Then the counit of the adjunction $\mathcal{E}\to \delta\vert \mathcal{E}\vert$ forgets to a weak equivalence, and hence is itself a weak equivalence. This proves the second claim.
\endofproof

To relate sheaves of chain complexes and simplicial presheaves, we will use the Dold--Kan correspondence. Let $\mathscr{C}{\rm h}_{\geq 0}$ denote the category of nonnegatively graded chain complexes. Recall the Dold--Kan functor $K:\mathscr{C}{\rm h}_{\geq 0}\to \sset$ denote the functor $K$ that sends a chain complex $C$ to the simplicial set
$$K(C)_n=\bigoplus_{[n]\onto [k]}C_k, \quad d_i(c_f)_{e(fd^i)}=m(fd^i)^*c_f, \quad s_i(c_f)_{fs^i}=c_f,$$
where given $f:[n]\to [m]$, $f=m(f)e(f)$ is the epi-mono factorization \cite[Section III]{GoerssJardine}. Recall also the normalized Moore functor $N:\sset\to \mathscr{C}{\rm h}_{\geq 0}$ defined on a simplicial set $X$ by 
$$N(X)_n=\bigcap_{i=0}^{n-1}\ker(d_n)\subset \ZZ(X)_n, \quad d=(-1)^nd_n,$$
where $\ZZ$ denotes the free simplicial abelian group functor. The pair of functors defines an adjunction $N\dashv K$. The following refinement of the Dold--Kan correspondence was observed in \cite[Proposition 2.2.31]{Urs} 
\begin{proposition}\label{DoldKan}
The pair of functors $N:\sm\mathscr{S}\to \sm\mathscr{C}{\rm h}_{\geq 0}$ and $K:\sm\mathscr{C}{\rm h}_{\geq 0}\to \sm\mathscr{S}$, obtained by applying the functors $N$ and $K$ objectwise, induces a Quillen adjunction $N\dashv K$. Moreover, if $C:\sm\mathscr{S}\to \sm\mathscr{C}{\rm h}_{\geq 0} $ denotes the alternating face map complex
$$C(X)_n=\ZZ(X)_n, \quad d=\sum_{i=0}^n (-1)^id_i,$$
then the inclusion $N(X)\into C(X)$ is a quasi-isomorphism of sheaves of chain complexes.
\end{proposition}

\section{A lemma about derived mapping spaces}

In this appendix, we prove a lemma that allows us to compute the derived mapping space as a discrete simplicial set when mapping into a zero-truncated object.
\begin{lemma}\label{truncate}
Let $X\in \sm\set\subset  \sm\mathscr{S}$ be a sheaf of smooth sets regarded as a simplicially discrete presheaf of simplicial sets.   Then for any $Y\in \sm\mathscr{S}$, the derived mapping space $R\map(Y,X)$ is weakly equivalent to the discrete simplcial set on $\set(\tilde \pi_0(Y),X)$.
 \end{lemma}
 \proof 
 Objectwise, $X$ is the nerve of a discrete groupoid. Hence, $X$ is an objectwise Kan complex. Since $X$ is a sheaf, it is a local object in $\sm\mathscr{S}$, hence a fibrant object in the projective model structure. Therefore, the derived mapping simplicial set can be computed in the simplicial model structure as $\map(QY,X)$. Explicitly, this simplicial set is given by 
 $$\map(QY,X)_n=\set(QY\times \underline{\Delta}^n,X),$$
 where $\underline{\Delta}^n$ denotes the constant simplicial presheaf on the simplicial $n$-simplex. Let $N:\set\to \mathscr{S}$ denote the nerve functor that sends a set to the corresponding discrete simplicial set. By the adjunction $\pi_0\dashv N$, we identify the $n$-simplices as 
 $$\set(QY\times \underline{\Delta}^n,X)\cong \set(\pi_0(QY\times \underline{\Delta}^n),X)\cong \set(\pi_0(QY),X),$$
 where the right isomorphism is given by distributing $\pi_0$ over the product and using the fact that $n$-simplices are connected. The above isomorphisms are clearly natural in $\underline{\Delta}^n$. Hence, we have an isomorphism of simplicial sets
 $$\map(QY,X)\cong N\set(\pi_0(QY),X).$$
 Since $X$ is a local object, the sheafification map $\pi_0(QY)\to \tilde \pi_0(QY)$ induces an isomorphism 
 $$\set(\tilde \pi_0(QY),X)\overset{\cong}{\to} \set(\pi_0(QY),X).$$
Now since the cofibrant replacement $QY\to Y$ is a local weak equivalence in the {\v C}ech model structure, it is also a local weak equivalence in the Jardine model structure. Hence, we have an induced isomorphism on sheaves
$\tilde \pi_0(QY)\to \tilde \pi_0(Y).$
Since $X$ is a sheaf, this implies that we have an induced isomorphism 
$$
\set(\tilde \pi_0(QY),X)\cong \set(\tilde \pi_0(Y),X)
$$ 
Composing with the above isomorphisms proves the claim.
\endofproof


\begin{thebibliography}{1}
\bibitem{Atiyah} M. Atiyah, F. Hirzebruch, {\it Vector bundles and homogeneous spaces}, Proc. Sympos. Pure Math. {\bf 3} (1961), 7–38.
\bibitem{BB} C. B\"ar, C. Becker, {\it Differential Characters}, Lecture Notes in Math. {\bf 2112}, Springer, New York,
(2014).
\bibitem{Bott} R. Bott, {\it On the Chern--Weil homomorphism and the continuous cohomology of Lie-groups}, Adv. Math. {\bf 2} (1973), 289-303.
\bibitem{Bunke} U. Bunke {\it Differential cohomology}, arXiv:1208.3961
\bibitem{BNV} U. Bunke, T. Nikolaus, M. V\"olkl, {\it Differential cohomology theories as sheaves of spectra}, J. Homotopy Relat. Struct. {\bf 11} (2016), 1-66.
\bibitem{CS} J. Cheeger, J. Simons, {\it Differential characters and geometric invariants}, Geometry and Topology, Lecture Notes in Mathematics {\bf 1167}, 50-80, Springer (1985).
\bibitem{Dugger} D. Dugger, {\it Spectral enrichments of model categories}, Homol. Homotopy Appl. {\bf 8} (1), 1 - 30 (2006)
\bibitem{Dug} D. Dugger, {\it Universal homotopy theories}, Adv. Math. {\bf 164} (1), 144-176 (2001).
\bibitem{Elmendorf} A. Elmendorf, M. Mandell, {\it Rings, modules and algebras in infinite loop space theory}, Adv. Math. {\bf 205} (1) (2006), 163-228.
\bibitem{FH} D. Freed, M. Hopkins, {\it Chern–Weil forms and abstract homotopy theory}, Bull. Am. Math. Soc {\bf 50} (3) (2013), 431–468.
\bibitem{GoerssJardine} P. Goerss, J. Jardine, {\it Simplicial homotopy theory}, Progress in Mathematics, Birkh\"auser (1999), Modern Birkh\"auser Classics (2009).
\bibitem{GradySati} D. Grady, H. Sati {\it Differential KO-theory: constructions, computations, and applications}, Adv. Math. {\bf 384} (2021), 107671.
\bibitem{Greenlees} J. Greenlees, {\it Equivariant connective $K$-theory for compact Lie groups}, J. Pure Appl. Algebra {\bf 187} (2004), 129 – 152.
\bibitem{Hirschhorn} P. Hirschhorn,  {\it Model categories and their localizations}, Mathematical Surveys and Monographs {\bf 99}, (2003).
\bibitem{HS} M. Hopkins, I. Singer, {\it Quadratic Functions in Geometry, Topology, and M-Theory}, J. Differ. Geom. {\bf 70} (3) (2005), 329-452.
\bibitem{HTT} J. Lurie, {\it Higher topos theory (AM-170)}, Princeton University Press, (2009).
\bibitem{LurieHA} J. Lurie {\it Higher Algebra}, https://people.math.harvard.edu/~lurie/papers/HA.pdf
\bibitem{MMSS} M. Mandell, P. May, S. Schwede, B. Shipley, {\it Model categories of diagram spectra}, Proc. London Math. Soc, {\bf 82} (2001), 441-512.
\bibitem{Pavlov} D. Pavlov J. Scholbach, {\it Admissibility and rectification of colored symmetric operads}, J. Topol. {\bf 11} (3) (2018), 559-601.
\bibitem{Urs} U. Schreiber, {\it Differential cohomology in a cohesive $\infty$-topos}, arXiv:1310.7930.
\bibitem{ShipSchwede} S. Schwede, B. Shipley, {\it Stable model categories are categories of modules}, Topology {\bf 42} (2003), 103-153.
\bibitem{Schipley} B. Shipley, {\it H$\ZZ$-algebra spectra are differential graded algebras}, Am. J. Math. {\bf 129} (2007), 351-379.
\end{thebibliography}
\end{document}